\documentstyle[12pt]{article}
\title{Proof of generalized Riemann hypothesis for Dedekind zetas and
Dirichlet L-functions}
\author{Andrzej M\c{a}drecki\thanks{Institute of Mathematics,
Wroc{\l}aw University of Technology (WUT), 50-370 Wroc{\l}aw, Poland}}
\newtheorem{de}{Definition}
\newtheorem{lem}{Lemma}
\newtheorem{th}{Theorem}
\newtheorem{pr}{Proposition}

\newtheorem{re}{Remark}

\newfont{\lll}{msbm10 scaled 1095}
\def\LA{\mbox{\lll \char65}}
\def\LC{\mbox{\lll \char67}}
\def\LF{\mbox{\lll \char70}}
\def\LH{\mbox{\lll \char72}}
\def\LN{\mbox{\lll \char78}}
\def\LQ{\mbox{\lll \char81}}
\def\LR{\mbox{\lll \char82}}
\def\LS{\mbox{\lll \char83}}
\def\LT{\mbox{\lll \char84}}
\def\LZ{\mbox{\lll \char90}}
\begin{document}
\maketitle

{\bf Abstract}. A short proof of the generalized Riemann hypothesis
(gRH in short) for zeta functions $\zeta_{k}$ of algebraic number
fields $k$ - based on the Hecke's proof of the functional equation
for $\zeta_{k}$ and the method of the proof of the Riemann hypothesis
derived in [$M_{A}$] (algebraic proof of the Riemann hypothesis) is given.
The generalized Riemann hypothesis for Dirichlet L-functions is an
immediately consequence of (gRH) for $\zeta_{k}$ and suitable product
formula which connects the Dedekind zetas with L-functions.

\section{Introduction}
Let $k$ be an {\bf algebraic number field}, (i.e. the main half of the set
of {\bf global fields}), i.e. a finite  algebraic extension of the
{\bf rational number field} $\LQ$. Let $R_{k}$ be a ring of
{\bf algebraic integers} in $k$ ,i.e. a finitely-generated ring
extension - the integral closure - of the ring of {\bf integers} $\LZ$. Then,
the {\bf Dedekind zeta function} $\zeta_{k}$ for $k$ is well locally defined (cf.e.g.
[K, Chapter 7], [L,VIII.2] and [N, VII])  as the {\bf Dirichlet series}
\begin{equation}
\zeta_{k}(s)\;:=\;\sum_{0\ne I\in {\cal I}_{k}}\frac{1}{N(I)^{s}}\;,\;Re(s)>1,
\end{equation}
where by $\LC$ we denote the field of
all complex numbers and by $Re(s)$ and $Im(s)$ the {\bf real} and {\bf
imaginary} part of a complex number $s$, respectively. We denote the
{\bf group of all fractional ideals} of the {\bf Dedekind ring} $R_{k}$
by ${\cal I}_{k}$ (cf.e.g. [N]) and finally $N(I)$ denotes the {\bf absolute norm}
of the ideal $I$, i.e. the number of elements in $R_{k}/I$.

We remark at once that we only use classical
Dirichlet-Dedekind-Hecke theory, from the heroic period of German
mathematics, to obtain an exciting result : a proof of the
{\bf generalized Riemann Hypothesis}( $gRH_{k}$ in short) for algebraic number
fields $k$. Hecke theory posseses such depth, that its
classical tools are sufficient to obtain $(gRH_{k})$. For example,
probably one of the most characteristic properties of the theory of classical number
theory is that, one may embed a number field in the Cartesian
product of its completions at the {\bf archimedean points}, i.e. in a
Euclidean space. In more recent years (more precisely since Chevalley
introduced ideles in 1936, and Weil gave his adelic proof of the
Riemann-Roch theorem soon afterwards), it has been found most
convenient also to take the product over the {\bf non-archimedean points}, with
a suitable restriction on the components - the {\bf adele ring}
$\LA _{k}$. However, we do not use the adele techique of {\bf Tate's
thesis} in this paper but stress {\bf Hecke's theory} and we do
not use the new achievements of algebraic number theory connected with
adeles and ideles.

When, we work with Dedekind zetas, it is surprising that at once we
obtain a very expanded apparatus of notions of the queen of mathematics -
algebraic number theory.

The main property of $\zeta_{k}$ is the existence of the following {\bf
Hecke - Riemann analytic continuation functional equation} (HRace in short,
cf.e.g. [L,XIII.3, Th.3])
\begin{displaymath}
\zeta_{k}^{*}(s)\;:=\;\frac{\mid d(k)\mid^{s/2}}{2^{r_{2}s}\pi^{ns/2}}
\Gamma(\frac{s}{2})^{r_{1}}\Gamma(s)^{r_{2}}
\zeta_{k}(s)\;=\;\frac{2^{r_{1}}h(k)R(k)}{w(k)s(s-1)} \;+\;
\end{displaymath}
\begin{displaymath}
\;+\;\sum _{0\ne I \in {\cal I}_{k}}\int_{\mid\mid y \mid\mid \ge 1}exp(-\pi
d(k)^{-1/n}N(I)^{2/n}Tr(y))[\Pi(y)^{s/2}\;+\;\Pi(y)^{(1-s)/2}]\frac{dy}{y}
\end{displaymath}
, where : $d(k)$ is the {\bf discriminant} of a field $k$
{(cf.e.g.[N,II.2]),}

$\Gamma(s) := \int_{0}^{\infty}e^{-x}x^{s-1}dx\;\;Re(s)>0$ is the
(classical) {\bf gamma function},

$r_{1}$ is the {\bf number of real embeddings} of $k$ into $\LC$,

$r_{2}$ is {\bf half of the number of complex embeddings} of $k$ into $\LC$,
(The pair $r=[r_{1},r_{2}]$ is called the {\bf signature} of $k$).

$h(k)$ is the {\bf class number} , $R(k)$  is the {\bf regulator of $k$}
(cf.[N, III.2]) and $w(k)$ is the number of {\bf roots of unity} lying
in $k$.

$S_{\infty}(k)$ denotes the set of {\bf archimedean absolute
values} of $k$,
$ n = n(k) = [k:\LQ]$  is the {\bf degree} of $k$ over $\LQ$,

$N_{v}(k)= N_{v}$ is the {\bf local degree} of $k$, which is 1 if $v$ is
a real point of $k$ and 2 if $v$ is a complex valuation from the set
$S_{\infty}(k)$.

Finally

$Tr_{k}(y):=\sum_{v\in S_{\infty}(k)}N_{v}y_{v}$ and $\prod( y ):= \prod_{v\in S_{\infty}(k)} y_{v}^{N_{v}}$.

From the topological point of view the answer to the question : where
are zeros and poles of $\zeta_{k}$ located - the algebraic number
theory {\bf characteristics ( arithmetics invariant)} : $d(k), r_{1}, r_{2},
n(k), h(k), R(k), w(k), S_{\infty}(k)$ - which appears in (HRace) - (as we
will show below) - are not so important, apart from the {\bf topological
invariants of $k$} , the signature $r(k)$, degree $n(k)$ and polynomial
$s(s-1)$.
For example, the invariants $h(k)
, R(k),w(k)$ and $r_{1}$ appear when we consider the {\bf residue value}
of $\zeta_{k}$ at the pole $s= 1$, but not when we consider the location of
the {\bf single pole} $\{1\} = I(\LC)\cap R(\LC)$, where the algebraic varietes $I(\LC)
:= \{s=u+iv\in \LC: v(1-2u)=0\}$ and $R(\LC):= \{s = u+iv\in \LC :
u(u-1)-v^{2} =0\}$ {\bf do not even depend on} $k$. Moreover, for the
purposes of this paper it is only important that $h(k)$ is {\bf
finite}, but the value of $h(k)$ is not itself important.
More exactly, we derive an essential generalization of (HRace), where
the $n$-dimensional {\bf  standard Gaussian function}
\begin{equation}
G_{n}(x)\;:=\; e^{-\pi \mid \mid x \mid \mid_{n}^{2}}\;,\; x \in \LR^{n}
\end{equation}
(here $\mid\mid .\mid\mid_{n}$ is the {\bf Euclidean norm} on $\LR^{n}$ and
obviously here, and all in the sequel, $\LR$ stands for the field of real
numbers), will be replaced by any smooth fixed point of ${\cal F}_{n}$.

The function $G_{n}$ is a {\bf fixed point} of the {\bf Fourier transform}
${\cal F}_{n}$ on the {\bf Schwartz space} ${\cal S}(\LR^{n})$ of smooth and
rapidly decreasing functions. If we replace $G_{n}$ by any other fixed point
$\omega_{+}$ of ${\cal F}_{n}$ from ${\cal S}(\LR)$, then we can extend the
(HRace) to the {\bf Fixed point Hecke Riemann analytic continuation equation}
(Face in short) (cf. Section 2).

The idea of the generalization of (HRace) to (Face) is, in some small sense
very similar to {\bf Grothendieck's} magnificent idea of the
generalization of the notion of set theory topology to category
topologies (e.g. the well-known {\bf etale cohomologies}) - to obtain
the required results : to prove $(gRH_{k})$ in our case and to prove the
{\bf Riemann hypothesis} for {\bf congruence Weil zetas}, respectively.

The following very important {\bf rational function} ( the {\bf polar-zero
part}) appears in HRace.
\begin{equation}
\;\;W_{k}(s)\;:=\;\frac{\lambda_{k}}{s(s-1)}(\;:=\;\frac{2^{r_{1}}h(k)R(k)}
{w(k)s(s-1)})\;;\;s\in \LC.
\end{equation}
Hence, in $W_{k}(s)$ is written a very important polynomial  $I$ of
two variables,  which does not depend on $k$! , with coefficients in $\LZ$:
\begin{equation}
\;\;I(s)\;:=\;Im(W_{k}(s))\mid s(s-1)\mid^{2}/\lambda_{k}\;=\;v(2u-1);\;
s=u+iv,\;u,v \in \LR.
\end{equation}
The function $I(s)$ is mainly responsible for the form of the {\bf
generalized Riemann Hypothesis} for $\zeta_{k}$ ($(gRH_{k})$ in short), i.e.
the following well-known implication ( as in the case of the Riemann
hypothesis , cf.[$M_{A}$]):
\begin{displaymath}
(gRH_{k}) \;\;If\;\zeta_{k}(s)=0\;and\;Im(s)\ne 0\;, then\;Re(s)=1/2 .
\end{displaymath}

According to (1.4) , the following {\bf Trivial Riemann Hypothesis} ((TRH)
in short) holds:
\begin{equation}
(TRH)\; If\; I(s) = 0\; and\; Im(s) \ne 0\;, then \;Re(s) =1/2.
\end{equation}

As in [$M_{A}$] we pose the following {\bf Algebraic conjecture
for $\zeta_{k}$}:
\begin{displaymath}
(TRH)\; implies\; (gRH_{k}).
\end{displaymath}

More exactly, let us consider the {\bf algebraic $\LR$-variete} $I(\LC)
:= \{s\in \LC : I(s) = 0\}$ and the zero-dimensional {\bf holomorphic
manifold} $\zeta_{k}(\LC) := \{s \in \LC : \zeta_{k}(s) = 0\}$. Then
the Riemann hypothesis $(gRH_{k})$ is a kind of relation between the
cycles (of $\LR^{2}$ and $\LC$, respectively) :  $I(\LC)$ (which does
{\bf not depend} on $k$) and $\zeta_{k}(\LC)$, i.e.
\begin{displaymath}
\;\;\;\zeta_{k}(\LC)\;\subset \;I(\LC).
\end{displaymath}

In the sequel, the bi-affine-linear form $I(u,v)$ of two real variables,
we call the {\bf fundamental form} of the class
$\{\zeta_{k} : k\:\;is\;an\;algebraic\;number\;field \}$.

Thus, {\bf topological information} on the isolated points of the
meromorphic function $\zeta_{k}$ is written - in fact - in the
algebraic varieties $I(\LC)$ and $I(\LC)\cap R(\LC)$, and therefore
there exists some {\bf unexpected} (and hence {\bf deep}) relation
between the {\bf arithmetic} of $I \in \LZ[u,v]$ over $\LR$ and the
{\bf arithmetic} of $\zeta_{k}$ over $\LC$. Moreover, the "serious"
$(gRH_{k})$ could be reduced to the formal consequence of the "non-serious"
(TRH) by calculating different kinds of integrals ( with respect to
different {\bf Haar measures}), which leads to the {\bf
subsequence} functional equation : let $Gal(\LC/\LR) = \{id_{\LC}, c\}$ be
the {\bf Galois
group} of $\LC$, i.e. $id_{\LC}$ denotes the identity automorphism of
$\LC$ and $c$ is the {\bf complex conjugation} automorphism :
\begin{equation}
\;\;\;c(z)\;=\;c(u+iv)\;:=\;u\;-\;iv ,
\end{equation}
which is an {\bf idempotent map}, i.e. $c^{2} = id_{\LC}$.
The following {\bf generalized Riemann hypothesis functional equation}
($(gRhfe_{k})$ in short) with a {\bf rational term} $I$ and the {\bf
action} of $Gal(\LC/\LR)$ indicates some "hidden" Galois symmetry of
$\zeta_{k}$ :
\begin{displaymath}
(gRhfe_{k})\;\;\;Im(\sum_{g \in
Gal(\LC/\LR)}(F_{g}\zeta_{k})(g(s)))\;=\;\frac{\lambda_{k}(f_{1}(s)
-f_{2}(s))I(s)}{\mid s(s-1) \mid^{2}} \;,\;Re(s)\in [0,1/2).
\end{displaymath}
In opposite to the $(gRhfe_{k})$ , the (HRace) gives an "open symmetry" of
$\zeta_{k}^{*}$ :
\begin{equation}
\zeta_{k}^{*}(s)\;=\;\zeta_{k}^{*}(1-s).
\end{equation}

As in the case of the Riemann hypothesis, the functional equation
$gRhfe_{k}$ - immediately implies the generalized Riemann hypothesis for the
Dedekind zetas due to TRH.
In comparison to $[M_{A}]$, we have significantly shorted the technicality of
the proof of the theorem on existence of $n$-dimensional {\bf RH-fixed points}.
We consider the non-commutative field of {\bf quaternions} $\LH$, endowed with
the Hilbert transform ${\cal H}_{\LH}$ of a measure $\mu$ (see Sect.3)
\begin{displaymath}
({\cal H}_{\LH}\mu)(h)\;:=\; \int_{\LH}\frac{d\mu (x)}{\mid \mid h-x
\mid\mid_{4}^{4}},
\end{displaymath}
and the product ring (with zero divisors) $\LQ_{p}\times \LQ_{q}$ of
different p-adic number fields endowed with the Hilbert transform
${\cal H}_{pq}$ :
\begin{displaymath}
({\cal H}_{pq}\mu)(a)\;:=\;\int_{\LQ_{p}\times \LQ_{q}}\frac{d
\mu(x)}{\Delta_{pq}(a-x)}.
\end{displaymath}

Thus, using the techniques used in [$M_{A}$] for the proof of the Riemann
hypothesis, we show that our method initiated in that article works and can
be significantly extended to the general case : this technique of RH-fixed
points - leads to the proof of the generalized Riemann hypothesis for Dedekind
zetas and Dirichlet L-functions.

The constructions in Section 3 are much more abstract in comparing to
$[M_{A}]$ and much simpler. Moreover, these construct are interesting in
themselves, since they  (and in some sense return) to fundamental
problems raused at the beginning of the 20th century.

The "heart" of the proof of RH from $[M_{A}]$ moving (practically
without any changes) for $gRH_{k}$.

\section{Fixed point Hecke-Riemann functional continuation equations}

These two chapters achieve two goals simultaneously. We present here all the
necessary preliminaries and notation. Next, we state the extension of (HRace)
to (Face). Secondly, the main technical tool - and in fact - the "heart of
the paper" , is Theorem 2 on the existence of multidimensional RH-fixed
points. Moreover, we comnent on a surprising property of the construction
mentioned: that it violates the Tertium non Datur in the
case, when the {\bf amplitude} $A$ has a {\bf support outside a set
of Lebesgue measure zero}.

Let $n \in \LN^{*}:=\LN-\{0\}$ be arbitrary (in all the sequel $\LN^{*}$
denotes the set of all positive integers). In the sequel $n=n(k)$ will
always be considered as the degree of a fixed algebraic number field
$k$, i.e. $n = [k:\LQ]$.

Exactly $n$ different embeddings of $k$ into the complex field $\LC$
exists. Indeed, by {\bf Abel's theorem} $k$ can be written in the
form $k = \LQ(a)$ for a suitable {\bf algebraic} $a$.

If $a_{1}, ... ,a_{n}$ are all complex roots of the {\bf minimal
polynomial} for $a$ over $\LZ$, then the mappings $C_{j}, j=1, ... ,n$ (
the {\bf conjugates of $k$}) defined by
\begin{equation}
C_{j}(\sum_{k=0}^{n-1}A_{k}a^{k})\;:=\;\sum_{k=0}^{n-1}A_{k}a_{j}^{k}
\end{equation}
(for $A_{0}, ... , A_{n-1} \in \LQ$) are all isomorphisms of $k$ {\bf
into} $\LC$, and every such isomorphism has to be of this form.

The fields $C_{j}(k)$ are called the {\bf fields conjugated} with $k$.

If $C_{j}(k) \subset \LR$, then it is called a {\bf real embedding} and
otherwise $C_{j}(k)$ is called a {\bf complex embedding}.

Note that if $C_{j}$ is {\bf complex}, then $c \circ C_{j}$ is again an
embedding, complex of course, and so the number of complex embeddings
is {\bf even}.
The number of such pairs of embeddings is usually denoted by $r_{2}(k) =r_{2}$
, and the number of {\bf real embeddings} by $r_{1}(k) = r_{1}$.

The pair $r = r(k) =[r_{1}, r_{2}]$ is called the {\bf signature} of $k$
(cf.e.g. [N, II.1])

We denote the {\bf Lebesgue measure} on $\LR^{n}$ , and the {\bf Lebesgue
measure} of $\LC^{n}$ by $d^{n}x$ and $d^{n}z$, respectively.

If $r=r(k) =[r_{1},r_{2}]$ is the {\bf signature of $k$}, then we
define the {\bf signature group} $G_{r}$ of $k$ as the product
\begin{equation}
G_{r}:=\LR_{+}^{r_{1}}\times(\LC^{*})^{r_{2}},
\end{equation}
of $r_{1}$ - exemplars of the multiplicative group $\LR^{*}_{+}$ of
{\bf positive real numbers} and $r_{2}$-exemplars of the {\bf
multiplicative group of complex numbers} $\LC^{*}$.

Obviously, $G_{r}$ is a Locally Compact Abelian group (LCA in short).
Hence, the {\bf Haar measure} is well defined. Its {\bf standardly normalized
Haar measure} will be denoted by $H_{r}$. It is well-known that $H_{r}$
is the product of the form :
\begin{equation}
dH_{r}(g)=\frac{d^{r_{1}}x}{\mid x \mid}\otimes
\frac{d^{r_{2}}z}{\mid z
\mid^{2}}=\otimes_{i=1}^{r_{1}}\frac{dx_{i}}{x_{i}}\otimes_{j=1}^{r_{2}}
\frac{dz_{j}}{\mid z_{j} \mid^{2}}.
\end{equation}
The signature group $G_{r}$ is obviously the multiplicative subgroup of the
{\bf Euclidean ring}
\begin{equation}
E_{r}\;:=\;\LR^{r_{1}}\times \LC^{r_{2}}\;\simeq \LR^{n},
\end{equation}
with the componentwise multiplication. It is obviously a ring with divisors
of zero.
In particular, $E_{r}$ has got the {\bf Haar module}  $\Delta_{r} = mod_{r}$
with the property
\begin{equation}
\Delta_{r}(g)\;=\;mod_{r}(g)\;=\;\prod_{i=1}^{r_{1}}\mid x_{i}\mid
\prod_{j=1}^{r_{2}}\mid z_{j} \mid^{2} \;,\;g=(x_{1}, ...
,x_{r_{1}},z_{1}, ... ,z_{r_{2}}).
\end{equation}
is well defined on $E_{r}$.
Moreover
\begin{equation}
dH_{r}(g)\;=\;\frac{d^{r_{1}}x \otimes d^{r_{2}}z}{mod_{r}(g)}.
\end{equation}
We denote the {\bf $mod_{r}$-unit sphere} of $G_{r}$ by $G_{r}^{0}$, i.e.
\begin{equation}
G_{r}^{0}\;:=\;\{g \in G_{r} : mod_{r}(g) =1\}.
\end{equation}
It is an elementary fact that we can write $G_{r}$ as the product
\begin{equation}
G_{r}\;=\;\LR^{*}_{+} \times G_{r}^{0},
\end{equation}
because any $g \in G_{r}$ can be written uniquely as
\begin{equation}
g\;=\;t^{1/n}c
\end{equation}
with $t \in \LR^{*}_{+}$ and $c \in G_{r}^{0}$. Here $c = \{c_{v}\}$
and $t^{1/n}c := (mod_{r}(c)^{1/n}\cdot (\frac{c_{v}}{mod_{r}c}))$.

We denote the {\bf Haar measures} of $G_{r}^{0}$ by $H_{r}^{0}$.
According to (2.15), the Haar measure $H_{r}$ can be considered as the
product of the {\bf Lebesgue measure} $dt/t$ on $\LR^{*}_{+}$ and the
appropriate {\bf Haar measure} $H_{r}^{0}$ on $G_{r}^{0}$.

For a large class of {\bf $\Gamma_{r}$-admissible} functions
$f:G_{r}\longrightarrow \LC$ the ($n$-dimensional) {\bf Mellin
transform} $M_{n}(f)$ or rather the {\bf signature Gamma} $\Gamma_{r}(f)$
(associated with $f$) is well-defined as
\begin{equation}
\Gamma_{r}(f)(s)\;:=\;\int_{G_{r}}mod_{r}^{s}(g)f(g)dH_{r}(g)
\;=:\;M_{n}(f)(s)\;,\;Re(s)>0.
\end{equation}

Recall that $f:\LR^{n}\longrightarrow \LC$ belongs to the {\bf Schwartz
space} ${\cal S}(\LR^{n})$ of rapidly decreasing functions, if for each
$n$-tuple of integers $\ge 0, k=(k_{1}, ... ,k_{n})$ and $l = (l_{1},
... , l_{n})$
\begin{displaymath}
\;\;\;p_{k,l}(f)\;:=\;sup_{x\in \LR} \mid x^{k}(D^{l}f)(x)\mid
<+\infty,
\end{displaymath}
where  $x^{k} := x_{1}^{k_{1}}...x_{n}^{k_{n}}$  and $D^{l} :=
D_{1}^{l_{1}}... D_{n}^{l_{n}}$, is a partial differential
operator.

It is easy to check (cf.e.g. [$M_{A}]$, Sect.2, Lemma1]) that the
following holds for $f \in {\cal S}(\LR^{n})$ :
\begin{equation}
\;\;\;\Gamma_{n}(f)(s)\;\in \;\LC\;\;if\;\;Re(s)>0,
\end{equation}
since ${\cal S}(\LR)\otimes ... \otimes {\cal S}(\LR)$ (n-times), is
dense in ${\cal S}(\LR^{n})$.

We denote the ($n$-dimensional) {\bf Fourier transform} of $f$ by
${\cal F}_{n}f$ (for ${\cal F}$-admissible functions):
\begin{equation}
\;\;{\cal F}_{n}(f)(x)\;:=\;\int_{\LR^{n}}e^{2\pi
ixy}f(y)d^{n}y\;=:\;\hat{f}(x)\;;\;x \in \LR^{n} ,
\end{equation}
where $xy := \sum_{k=1}^{n}x_{i}y_{i}$ is the standard euclidean scalar product
of $n$-vectors $x = (x_{1}, ... ,x_{n})$ and $y = (y_{1},... ,y_{n})$.
In this paper, it is also very convenient to use the (1-dimensional)
{\bf plus-Sin transform} defined as
\begin{equation}
\;\;\;S_{+}(f)(x)\;:=\;\int_{0}^{+\infty}sin(xy)f(y)dy
\;=:\;\hat{f}_{+}(x)\;:\;x \in \LR_{+}.
\end{equation}

For another large class of {\bf $\theta$-admissible} functions
$f: G_{r} \longrightarrow \LC$ ($n$-dimensional or signatural), the{\bf Jacobi
theta function $\theta_{r}(f)$ associated with $f$} is defined as the series
\begin{equation}
\theta_{r}(f)(x)\;:=\;\sum_{k \in (\LN^{*})^{n}}f(k \cdot
x)\;=\;\int_{(\LN^{*})^{n}} f(k \cdot x)dc(x)\;,x\in \LR_{+}^{n},
\end{equation}
where $k \cdot x$ denotes componentwise multiplication in $E_{r}$
and $dc$ is the {\bf calculating measure} on $(\LN^{*})^{n}$ , i.e. the
unique Haar measure on $\LZ^{n}$ normalized by the condition :
$c(\{0\})=1$.

Beside the field $\LC$, we will also use the non-commutative field of
quaternions $\LH$. It is well-known (cf.e.g. [W]) that the formula
\begin{equation}
\Delta_{\LH}(h)\;:=\;\mid\mid h \mid\mid_{4}^{4} \;,h \in \LH ,
\end{equation}
defines the {\bf Haar module} of $\LH$.

For a class of some ${\cal H}$-admissible measures defined on a compact
subset $C$ of $\LH$, we define the (compact) $\LH$-{\bf Hilbert transform}
${\cal H}_{\LH}$ by the formula
\begin{equation}
({\cal H}_{\LH}\mu)(h)\;:=\;\int_{C}\frac{d\mu(x)}{\Delta_{\LH}(h-x)}\;,\;h
\in \LH.
\end{equation}
Finally, we use the product ring $\LQ_{p}\times \LQ_{q}$ with zero
divisors of different p-adic number fields. It is well-known that the
formula
\begin{displaymath}
\Delta_{pq}(x_{p},x_{q})\;:=\;\mid x_{p} \mid_{p} \mid x_{q}
\mid_{q}\;,(x_{p},x_{q})\in \LQ_{pq},
\end{displaymath}
defines the {\bf Haar module} of $\LQ_{pq}$ and the formula
\begin{displaymath}
({\cal H}_{pq}\mu)(a)\;:=\;\int_{\LQ_{p}\times \LQ_{q}} \frac{d
\mu(x)}{\Delta_{pq}(a\;-\;x)}\;,\; a \in \LQ_{pq},
\end{displaymath}
defines pq-Hilbert transform.

Finally, we note that the Schwartz spaces ${\cal S}(\LR^{n})$ are
{\bf admissible} for all the integral transforms defined above :
$\Gamma_{r}, {\cal F}_{n}, \theta_{r}$ and ${\cal H}$ ( for absolutely
continuous measures $\mu$ w.r.t. Lebesgue measure $d^{4}h$ and the Haar
measure $dH_{pq}$ of $\LQ_{pq}$, considered as densities of signed measures).

One of the main tools when we work with zetas is the {\bf Poisson
Summation Formula} (PSF in short, cf.e.g. [N], [L, XIII.2]) , which
shows that ${\cal F}_{n}$ is a $l^{1}(\LZ)$-{\bf quasi-isometry} on
${\cal S}(\LR^{n})$ and using our notation can be written as :
\begin{displaymath}
(PSF)\;\;\;\int_{\LZ^{n}}\hat{f}(x)dc(x)\;=\;\int_{\LZ^{n}}
f(x)dc(x),
\end{displaymath}
if $f \in {\cal S}(\LR^{n})$.

A complex function $\omega_{+}$ on $\LR^{n}$ ($n = r_{1}+2r_{2}$) is called
a {\bf fixed point of ${\cal F}_{n}$}, if it is an {\bf eigenvector} of
${\cal F}_{n}$ with the corresponding {\bf eigenvalue} equal to 1, i.e.
\begin{equation}
{\cal F}_{n}(\omega_{+})\;=\;\hat{\omega_{+}}\;=\;\omega_{+}.
\end{equation}

Analogously a complex valued function $\omega_{-}$ on $\LR^{n}$ is called
the {\bf -fixed point} of ${\cal F}_{n}$ if it is an {\bf eigenvector} of
${\cal F}_{n}$ corresponding to the {\bf eigenvalue} $-1$ of ${\cal
F}_{n}$ :
\begin{displaymath}
{\cal F}_{n}(\omega_{-})=\hat{\omega_{-}}=-\omega_{-}.
\end{displaymath}

We use the common name for $+$fixed points and $-$fixed point - the
{\bf $\pm$fixed points}  $\omega_{\pm}$.

Let $\omega = \omega_{\pm}$ be a {\bf $\pm$fixed point} of ${\cal F}_{n}$
from ${\cal S}(\LR^{n})$ and let $M =[m_{ij}]_{n\times n}$ be a {\bf matrix}
of real numbers.

Let us consider the function
\begin{displaymath}
\omega_{M}(x)\;:=\;\omega_{\pm}(Mx^{t})\;\;;\;\;x^{t} \in \LR^{n},
\end{displaymath}
and the theta associated with it
\begin{equation}
\theta_{n}(\omega_{M})(x)\;:=\;\sum_{m\in \LZ^{n}} \omega_{M}(mx)\;,\;x \in
\LR^{n}.
\end{equation}

\begin{lem}({\bf Hecke's theta formula})

For each non-singular matrix $M$ the following relation holds
\begin{displaymath}
(HTF)\;\;\;\theta_{n}(\omega_{M}^{\pm})(x)\;\;=\;\;\pm \theta_{n}
(\omega_{^{t}M^{-1}}^{\pm})(x)/\mid det(M) \mid.
\end{displaymath}
\end{lem}
{\bf Proof}. Let $M = [m_{ij}]_{n\times n}$ and $\omega_{M}(x) :=
\omega(Mx^{t}), x \in \LR^{n}$. If $M$ is a non-singular real matrix, then
$det(M)\ne 0$. Thus the function $\omega_{M}$ is also in ${\cal
S}(\LR^{n})$ and using the change of variables formula for multiple
integrals, we immediately find that its Fourier transform is
given by
\begin{displaymath}
\hat{\omega_{M}^{\pm}}(x)\;=\;\pm \frac{\omega(^{t}M^{-1}x^{t})}
{\mid det(M) \mid},
\end{displaymath}
where $^{t}M^{-1}$ is the transpose of the inverse of $M$.

This is clear, since when we make the change of variables $z = Mx^{t}$, we
have $dz =\mid det(M) \mid dx$, and  $<M^{-1}z^{t}, y> = <z, ^{t}M^{-1}y^{t}>$.

The first important step in the proof of $(gRH_{k})$ is
the {\bf generalization of the Hecke-Riemann analytic continuation
eqation} ((HRace) in short), given below. Therefore we need some additional
notation.

Let us again consider the signature {\bf Euclidean ring}
\begin{displaymath}
E_{r}\;=\;\LR^{r_{1}}\times \LC^{r_{2}}\;\simeq \LR^{n},
\end{displaymath}
and the {\bf conjugation} map $C : k \longrightarrow E_{r}$ defined as
\begin{equation}
C(\xi)\;:=\;(C_{1}(\xi), ... , C_{n}(\xi))\;;\xi \in k.
\end{equation}
Let us observe that each conjugate $C_{v}$ determines the {\bf absolute
value} (place) $v$ of $k$ by the formula :
\begin{equation}
v(\xi)\;:=\;\mid C_{v}(\xi) \mid\;;\;\xi \in k.
\end{equation}

The {\bf completion} of $(k,v)$ is denoted by $k_{v}$. Since $v$ is
{\bf archimedean}, $k_{v}$ is equal to $\LR$ or $\LC$.

In the case : $k_{v} \simeq \LC$ the completion is determined up to
{\bf complex conjugation} $c$, according to the well-known elementary
fact that if $\sigma \in Gal(\LC/\LR)$, then
\begin{equation}
v(\sigma(z))\;=\;v(z)\;,\;z \in \LC
\end{equation}
(cf.e.g. [L, II.1] and [N, L.3.1]).

We denote the set of all {\bf non-equivalent archimedean places} of $k$
by $S_{\infty}(k)$. According to (2.28), it is obvious that
\begin{displaymath}
\mid S_{\infty}(k)\mid \;=\;r_{1}\;+\;r_{2}.
\end{displaymath}

Let us consider the map $\mid C \mid : k \longrightarrow \mid E_{r}
\mid\;:=\;\LR^{r_{1}+r_{2}}$ defined as
\begin{equation}
\mid C \mid (\xi)\;:=\;(\mid C_{v}(\xi) \mid :v \in S_{\infty}(k)).
\end{equation}
Recall that $G_{r} = (\LR^{*}_{+})^{r_{1}}\times (\LC^{*})^{r_{2}}$.
So, if we denote : $\mid G_{r}\mid := (\LR_{+}^{*})^{r_{1}+r_{2}}$, then
we have the decomposition
\begin{equation}
G_{r}\;\simeq \;\mid G_{r}\mid \times \LT^{r_{2}},
\end{equation}
where $\LT :=\{z \in \LC : \mid z \mid =1\}$ is the 1-dimensional
torus.

The kernel of $\mid C \mid$, i.e. $\mu(k) := ker(\mid C \mid)$ is the
{\bf group of the roots of unity} in $k$. Let
\begin{equation}
w(k)\;=\;\# \mu(k)\;=\;\mid \mu(k) \mid,
\end{equation}
be the {\bf number of roots of unity in} $k$.

Let $U(k)$ be the {\bf group of units of $k$} ($S_{\infty}(k)-{\bf
units})$), i.e.
\begin{displaymath}
U(k)\;=\;R_{k}^{*} .
\end{displaymath}

Let $V(k) \;:=\;\mid C \mid(U(k))$ be the {\bf image} of $U(k)$ under
the mapping $\mid C \mid$. Its image $V(k)$ is contained in the subgroup
$\mid G_{r}^{0} \mid$ consisting of all $g \in \mid G_{r} \mid$ such that
$mod_{r}(g) = 1$, and is a {\bf discrete subgroup}. Furthermore, $\mid
G_{r}^{0} \mid/V(k)$ is {\bf compact} (cf. [L,p.256]). Also, we can
write $G_{r}$ as the product
\begin{displaymath}
G_{r}\;=\;\LR_{+}^{*}\times \mid G_{r}^{0} \mid \times \LT^{r_{2}}.
\end{displaymath}
Finally, let $E(k)$ be the {\bf fundamental domain for $V^{2}(k)$} in
$\mid G_{r}^{0} \mid$ (cf. [L]). We obtain the following {\bf disjoint
decomposition}
\begin{equation}
\mid G_{r}^{0} \mid \;=\;\cup_{\eta \in V} \eta^{2} E(k).
\end{equation}

Let $A$ be an arbitrary {\bf integral (fractional) ideal} of $k$. Then
, it is well-known that $A$ has an {\bf integral basis} over $\LZ$
(cf.[N, Th.2.4]). Thus, each $\xi \in A$ can be written as
\begin{equation}
\xi\;=\;x_{1}\alpha_{1}\;+\;...\;+\;x_{n}\alpha_{n}\;\;,\;\;x_{i}\in
\LZ.
\end{equation}

For $v \in S_{\infty}(k)$ we let $C_{v}$ be the embedding (conjugate)
of $k$ in $k_{v}$, identified with $\LR$ or $\LC$ (in the case of $\LC$,
we fix one identification, which otherwise is determined only up to
conjugacy). We will write
\begin{displaymath}
\xi_{v}\;=\;C_{v}(\xi)\;=\;\sum_{j=1}^{n}x_{j}C_{v}(\alpha_{j})
\end{displaymath}
and
\begin{displaymath}
C(A)\;:=\;[N(A)^{-1/n}C_{i}(\alpha_{j})] \;,i,j=1,...,n.
\end{displaymath}
Hence, $N(A)^{-1/n}[\xi_{1}, ... ,\xi_{n}] = C(A)[x_{1}, ... , x_{n}]^{t}$
and we also use  this same notation when we constrict $x_{i}$ to the set of
real numbers.

Let ${\cal R}$ be an class of ideals of the ordinary ideal class group $H(k)
:={\cal I}_{k}/P_{k}$. Let $A$ be an ideal in ${\cal R}^{-1}$. The map
\begin{equation}
B \longrightarrow  AB\;=\;(\xi)
\end{equation}
eatablishes a bijection between the set of ideals in ${\cal R}$ and
equivalence classes of non-zero elements of $A$ : $A/\sim_{u}$, where two
field elements are called {\bf equivalent} $\sim_{u}$, if they differ by a
{\bf unit}.

Let $R(A)$ be a set of {\bf representatives} for the non-zero
equivalence classes.

Finally, we introduce two thetas - small and capital : the {\bf small Jacobi
theta of $k$} (associated with $\omega$)
\begin{equation}
\theta_{k}(\omega)(g)\;\;:=\;\;\sum_{0\ne I \in {\cal I}_{k}}\sum_{\xi
\in R(I)}\sum_{u \in U(k)}\omega(u \xi g) \;=\;
\end{equation}
\begin{displaymath}
\;=\;\sum_{0\ne I \in {\cal I}_{k}} \sum_{x \in \LZ^{n}}
\theta_{n}(\omega_{C(I)})(g)\;;\;g\in G_{r},
\end{displaymath}
and the {\bf radial Jacobi theta} of $k$
\begin{displaymath}
\Theta_{k}(\omega)(t)\;:=\;\frac{\int_{E(k)}\theta_{k}(\omega)(ct^{1/n})dH_{r}
^{0}(c)}{w(k)}\;,\;t \in \LR_{+}^{*}.
\end{displaymath}

\begin {th}({\bf Fixed point HRace = Face})
The following functional equation holds for each $\pm$fixed point
$\omega_{\pm}$ of ${\cal F}_{n}$ from ${\cal S}(\LR^{n})$, with the property that
$\Gamma_{r}(\omega_{\pm})$ {\bf does not vanishes}, and for each $s$ with $Re(s)>0$
\begin{equation}
(Face)\;\;(\Gamma_{r}(\omega_{\pm})\zeta_{k})(s)\;=\;\frac{\lambda_{k}\ne 0}
{s(s-1)}\;+\;
\end{equation}
\begin{displaymath}
\;+\;\int_{1}^{\infty}\int_{E(k)}\theta_{k}(\omega_{\pm})(ct
^{1/n})(t^{s}\;\pm\;t^{1-s})dH_{r}^{0}(c)\frac{dt}{t}=\int_{1}^{\infty}
\Theta_{k}(\omega_{\pm}(t)(t^{s-1}+t^{-s}))dt.
\end{displaymath}
\end{th}
{\bf Proof}. ( A topological simplification of Lang's version of Hecke's proof
of (HRace)).
Let ${\cal R}$ be an {\bf ideal class} of the ordinary {\bf ideal class
group} $H(k) :={\cal I}_{k}/P_{k}$ , where $P_{k}$ is the subgroup of
principal fractional ideals.
It is convenient to deal at initially with the zeta function associated
with an ideal class ${\cal R}$.  We define
\begin{equation}
\zeta_{k}(s, {\cal R})\;:=\;\sum_{B \in {\cal R}} \frac{1}{N(B)^{s}}
\end{equation}
for $Re(s)>1$. Let $A$ be an ideal in ${\cal R}^{-1}$. Then the map
\begin{equation}
B \longrightarrow AB\;=\;(\xi)
\end{equation}
establishes a bijection between the set of ideals in ${\cal R}$ and
{\bf equivalence classes of non-zero elements of $A$} (where two field
elements are called equivalent, if they differ by a unit from $U(k)$).
Let $R(A)$ be a set of {\bf representatives} for the non-zero
equivalence classes. Then
\begin{equation}
N(A)^{-s}\zeta_{k}(s,{\cal R})\;=\;\sum_{\xi \in R(A)} mod_{r}(\xi
N(A)^{-1/n})^{-s}.
\end{equation}

We recall that the signature gamma is represented by the following integral
(cf.(2.17))
\begin{displaymath}
\Gamma_{r}(\omega_{\pm})(s)\;=\;\int_{G_{r}}\omega_{\pm}(g)
mod_{r}(g)^{s}dH_{r}(g),
\end{displaymath}
for $Re(s)>0$, since
\begin{displaymath}
mod_{r}(\xi)\;=\;\prod_{v \in S_{\infty}(k)}\mid \xi_{v} \mid^{N_{v}},
\end{displaymath}
where $N_{v} =[k_{v}:\LR]$ are local degrees.

It will also be useful to note that if $f$ is a function such that
$f(g)/ mod_{r}(g)$ is absolutely integrable on $G_{r}$, then
\begin{displaymath}
\int_{G_{r}}f(g)\frac{dH_{r}(g)}{mod_{r}(g)}\;=\;\int_{G_{r}}f(Mg)
\frac{dH_{r}(g)}{mod_{r}(g)},
\end{displaymath}
for any n-dimensional matrix $M = [m_{ij}]$ with {\bf real} $m_{ij}$.

In other words, $dH_{r}(g)/mod_{r}(g)$  is an {\bf invariant
measure} of the {\bf dynamical system} $(G_{r} , T_{M}(y) := My)$ or,
in other words, $H_{r}/\Delta_{r}$ is a {\bf Haar measure} on the group
$G_{r}$.

Note that the signatural gamma function is expressed as such an integral.

Therefore, substituting $g$ by $N(A)^{1/n}\xi g$ in (2.17), we obtain
\begin{equation}
\Gamma_{r}(\omega_{\pm})(s)\frac{N(A)^{s}}{mod_{r}(\xi)^{s}}\;=\;
\int_{G_{r}}\omega_{\pm}(\xi N(A)^{-1/n}g) mod_{r}(g)^{s}dH_{r}(g)
\end{equation}

For $Re(s)\ge 1+\delta$, the sum over inequivalent $\xi \ne 0$ is
absolutely and uniformly convergent. Since for $Re(s)>1$,
\begin{displaymath}
N(A)^{-s}\zeta_{k}(s, {\cal R})\;=\;\sum_{\xi \in R(A)}mod_{r}(\xi
N(A)^{-1/n})^{-s},
\end{displaymath}
it follows that
\begin{equation}
\Gamma_{r}(\omega)(s)\zeta_{k}(s,{\cal R})
=\int_{G_{r}}\sum_{\xi \in R(A)} \omega_{D(\xi)N(A)^{-1/n}}(g)
mod_{r}(g)^{s}dH_{r}(g),
\end{equation}
where $D(\xi) := [\delta_{iv}C_{v}(\xi)]$ denotes a diagonal matrix of
conjugations.
But according to (2.30),  we can write
\begin{displaymath}
g\;=\;t^{1/n}c \;,\;t>0,c \in G_{r}^{0}.
\end{displaymath}
Therefore,
\begin{equation}
\Gamma_{r}(\omega_{\pm})(s)\zeta_{k}(s,{\cal
R})=\int_{0}^{\infty}\int_{G_{r}^{0}}\sum_{\xi \in
R(A)}mod_{r}(t^{1/n}c)^{s} \omega_{\pm}(N(A)^{-1/n}(\xi_{1}t^{1/n}c_{1}, ...
,\xi_{n}t^{1/n}c_{n}))t^{s}dH_{r}^{0}(c)\frac{dt}{t},
\end{equation}
where $dH_{r}^{0}(c)$ is the appropriate measure on $G_{r}^{0}$ and $c
=(c_{v})$ is a variable in $G_{r}^{0}$.

According to the decomposition (2.30) and since the kernel of $\mid C
\mid$  is the group $\mu(k)$, we obtain from the above equation
\begin{equation}
\Gamma_{r}(\omega_{\pm})(s)\zeta_{k}(s,{\cal
R})=\int_{0}^{\infty}\int_{E(k)}\frac{t^{s}}{w(k)}\sum_{u \in U(k)}\sum_{\xi
\in R(A)}\omega_{\pm}(N(A)^{-1/n}(C_{1}(\xi u)t^{1/n}e_{1},...,C_{n}(\xi
u)t^{1/n}e_{n}))dH_{G_{0}}(e)dt/t=
\end{equation}
\begin{displaymath}
=\frac{1}{w(k)}\int_{0}^{\infty}\int_{E(k)}t^{s}\sum_{u \in U(k)}\sum_{x \in
X(A)}\omega_{\pm}((N(A)^{-1/n}(\sum_{j=1}^{n}x_{j}C_{1}(\alpha_{j}u)t^{1/n}e_{1}, ...
,\sum_{j=1}^{n}x_{j}C_{n}(\alpha_{j}u)t^{1/n}e_{n}))dH_{G_{0}}(e)dt/t ,
\end{displaymath}
where the second sum is over a {\bf subset} $X(A)$ of $\LZ^{n}-\{0\}$.
But, according to the definition of $R(A)$ , operating units we obtain that
if $u$ runs $U(k)$ and $\xi$ runs $R(A)$ then $x = [x_{1}, ... , x_{n}]
\in \LZ^{n}-\{0\}$  from
\begin{displaymath}
u\xi\;=\; \sum_{x \in \LZ^{n}-\{0\}} x_{j}\alpha_{j}
\end{displaymath}
spars all $\LZ^{n} - \{0\}$. Therefore, the "fourth integral" from
(2.43) we can rewrite in the form
\begin{equation}
=\int_{0}^{\infty}\int_{E(k)}\frac{t^{s}}{w(k)}\sum_{x=(x_{1}, ...
,x_{n})\ne 0}\omega_{\pm}((N(A)^{-1}t)^{1/n}e_{1}\sum_{j=1}^{n}x_{j}C_{1}
(\alpha_{j}, ... , (N(A)^{-1}t)^{1/n}e_{n}\sum_{j=1}^{n}C_{n}(\alpha_{j}))
)dH_{r}^{0}(e) dt/t=
\end{equation}
\begin{displaymath}
=\int_{0}^{\infty}\int_{E(k)}(\frac{t^{s}}{w(k)})\sum_{0\ne x \in
\LZ^{n}}\omega((N(A)^{-1}t)^{1/n}e C(A)x^{t})dH_{r}^{0}(e)dt/t \;=\;
\end{displaymath}
\begin{displaymath}
\;=\;\int_{0}^{\infty}\int_{E(k)}\frac{\theta_{n}(\omega_{C(A)})
(t^{1/n}e)-1}{w(k)}dH_{r}^{0}(e) dt/t .
\end{displaymath}

We split the integral from $0$ to $\infty$ into two integrals, from $0$
to $1$ and from $1$ to $\infty$. We thus find
\begin{equation}
\Gamma_{r}(\omega_{\pm})(s)\zeta_{k}(s,{\cal R})\;=\;\frac{1}{w(k)}
\int_{0}^{1}t^{s} \int_{E(k)}\theta_{n}(\omega_{C(A)}(t^{1/n}e))t^{s}dH_{r}^{0}
(e)dt/t\;-\;
\end{equation}
\begin{displaymath}
-\;\frac{H_{G_{r}^{0}}(E(k))}{w(k)s}\;+\;\int_{1}^{\infty}\int_{E(k)}
\frac{t^{s}}{w(k)}[\theta_{n}(\omega_{C(A)}(t^{1/n}e)-1]dH_{r}^{0}(e) dt/t .
\end{displaymath}
We return to the basis $\{\alpha_{j}: j=1, ... , n\}$ of the integral ideal
$A$ over $\LZ$. We define
\begin{displaymath}
\alpha^{*}\;:=\;\{\alpha_{j}^{*} : j =1, ... ,n\}
\end{displaymath}
to be the {\bf dual basis} with respect to the trace (cf. [L, XII.3]).
Then $\alpha^{*}$ is a basis for the fractional ideal
\begin{displaymath}
A^{*}\;:=\;(D_{k/\LQ} A)^{-1},
\end{displaymath}
where $D_{k/\LQ}$ is the {\bf different} of $k$ over $\LQ$ (cf. [L,
III.1], [N, IV.2] and the remark below).

We now use Heckes's theta functional equation $(HTE)$. It can be seen that
\begin{equation}
\theta_{n}(\omega_{C(A)}^{\pm})(t^{1/n}c)\;=\;\pm \frac{1}{t}\theta(\omega
^{\pm}_{C(A^{*})})(t^{-1/n}c^{-1}),
\end{equation}
because $mod_{r}(c) = 1$ , i.e. $c$ is in $G_{r}^{0}$ ! We transform the
first integral from $0$ to $1$, using a simple change of variables, letting
$t=1/\tau, dt = -d\tau/\tau^{2}$. Note that the measure $dH_{r}^{0}(c)$  is
invariant under the transformation $c \longrightarrow c^{-1}$ (think of
an isomorphism with the additive Euclidean measure, invariant under taking
negatives).

We therefore find that
\begin{equation}
\Gamma_{r}(\omega_{\pm})(s)\zeta_{k}(s,{\cal R})
\;=\;\frac{2H_{r}^{0}(E(k)\times  \LT^{r_{2}})}{w(k)s(s-1)}\;+\;
\end{equation}
\begin{displaymath}
\;+\;\frac{1}{w(k)}\int_{1}^{\infty}\int_{E(k)}(\theta_{n}(\omega_{C(A)}
)(t^{1/n}c)t^{s}\;\pm\;\theta_{n}(\omega_{C(A^{*})})(t^{1/n}c)t^{1-s})dH_{r}
^{0}(c)\frac{dt}{t}.
\end{displaymath}
(Let us remark that in the second edition of [L] in Section XII.3 , on page 257
there is a typegraphical error).

The expression in (2.47) is {\bf invariant} under the transformations

$A \longrightarrow A^{*}$ and $s \longrightarrow 1-s$ (in the plus case).

Thus, we have obtained full calculations on the zeta function of an ideal
class ${\cal R}$. Taking the sum over the ideal classes ${\cal R}$ from
$H(k)$ we immediately yield information on the zeta function itself, as
follows : we can construct for $A^{*}$  in a similar way ,and hence we finally
obtain
\begin{equation}
2(\Gamma_{r}(\omega_{\pm})\zeta_{k})(s)=\Gamma_{r}
(\omega_{\pm})(s)(\sum_{{\cal R}\in H(k)}2\zeta_{k}(s,{\cal R}))=
\end{equation}
\begin{displaymath}
=\frac{\lambda_{k}}{s(s-1)}+\frac{1}{w(k)}\int_{1}^{\infty}
\int_{E(k)}(t^{s}\pm t^{1-s}) (\sum_{{\cal R} \in H(k)}\sum_{A \in {\cal
R}^{-1}}\theta_{n}(\omega_{C(A)}(t^{1/n}c))d^{*}c \frac{dt}{t}=
\end{displaymath}
\begin{displaymath}
=\frac{\lambda_{k}}{s(s-1)}+\frac{1}{w(k)}\int_{1}^{\infty}\int_{E(k)}
(t^{s}\pm t^{1-s})\theta_{k}(\omega_{\pm})(t^{1/n}c)dH_{r}^{0}(c) \frac{dt}{t}.
\end{displaymath}
\begin{re}
As we mentioned above, in algebraic number theory we have to deal with
a very expanded notional aparatus. We recall some ideas, explored in
this paper.

Let $k$ be an arbitrary algebraic number field. Then we denote the {\bf trace}
of $k$ over $\LQ$ by $tr_{k}$.

If $A$ is a {\bf fractional ideal} of $k$, then $A^{*}$ denotes the
{\bf complementary ideal} to $A$ with respect to the trace $tr_{k}$,
defined as
\begin{equation}
A^{*}\;:=\;\{x \in k: tr_{k}(x A) \subset R_{k}\},
\end{equation}
(cf. [L, II.1]).
If $\{\alpha_{1}, ... , \alpha_{n}\}$ is a {\bf basis} of $A$ over
$\LZ$, then $\{\alpha_{1}^{*}, ... ,\alpha_{n}^{*}\}$
, where $\{\alpha_{i}^{*}\}$ is the dual basis relative to the
trace $tr_{k}$, is a basis of $A^{*}$.

One of the main notions of algebraic number theory is the {\bf different}
$D_{k/\LQ}$ . The different $D_{k/\LQ}$ "differs" $A^{-1}$ from
$A^{*}$, i.e. cf.e.g. [K], [L] and [N]
\begin{equation}
A^{*}\;=\;(D_{k/\LQ}A)^{-1}.
\end{equation}
One can show that
\begin{equation}
D_{k/\LQ}\;:=\;R_{k}^{*}.
\end{equation}
The second main important notion is the {\bf discriminant} $d(k)$ of an
algebraic number field.

If $\{C_{j}\}$ are embeddings as considered above and $\{\alpha_{j}\}$
forms a base of a fractional ideal $A$, then we can define the {\bf
discriminant}  $d_{k}(\alpha_{1}, ... , \alpha_{n})$  by
\begin{equation}
d_{k}(\alpha_{1}, ... ,
\alpha_{n})\;:=\;(det[C_{j}(\alpha_{i})]_{i,j})^{2}=det[tr_{k}(\alpha_{i}
\alpha_{j})].
\end{equation}
It is well-known that the discriminant of a basis of $A$ does not
depend on the choice of this basis. In particular, if $A = R_{k}$, then
this discriminant is called the {\bf discriminant of the field $k$} and
denoted by $d(k)$.
The discriminant $d(k)$ has many nice and important properties :

(1) according to the {\bf Stickelberger theorem}, $d(k)$ is either
congruent to unity (mod 4) or is divisible by 4,

(2) is strictly connected with the signature $r =[r_{1},r_{2}]$ :
$sign d(k) = (-1)^{r_{2}} $ , and according to the {\bf Minkowski
theorem}
\begin{displaymath}
\mid d(k) \mid\;>\;(\frac{\pi}{4})^{2r_{2}}(\frac{n^{n}}{n!})^{2} ,
\end{displaymath}
which also ilustates the strict relation with the degree $n =n(k)$.

(3) The connection with the different :
\begin{displaymath}
N(D_{k/\LQ})\;=\;\mid d(k) \mid.
\end{displaymath}
However the value of $d(k)$ is mainly underlined by the deep {\bf
Hermite theorem}, which asserts that only a finite number of algebraic fields
can have the same discriminant.

Besides the importance of $d(k)$, its arithmetic invariance does not appear
in our "topological" generalization of HRace.

We saw that one of the main roles in the proof of (Face) was played by
the function $\omega_{C(A)}$. In Lang's proof of HRace [L,XII.3],
this corresponds to the consideration of the gaussian fixed point $\omega
:= \otimes_{j=1}^{r_{1}}G_{1}\otimes_{j=1}^{r_{2}}G_{\LC}$, where
\begin{displaymath}
G_{\LC}(z)\;:=\;e^{-\pi \mid z \mid^{2}}\;;\;z \in \LC,
\end{displaymath}
is the {\bf complex Gaussian fixed point}  of ${\cal F}_{2}$ on $\LC
(=\LR^{2})$.

Then
\begin{displaymath}
\omega(C(A)x^{t})=exp(-\pi (N(A)^{2}d(k))^{-1/n}\sum_{j=1}^{n}\mid
\sum_{v=1}^{n}C_{v}(\alpha_{j})x_{j}\mid^{2})=:exp(-\pi(N(A)d(k))^{-1/n}
<A_{\alpha}x, x>),
\end{displaymath}
where the $\nu \mu$-component of the matrix $A_{\alpha} = [a_{\nu \mu}]$
is given by
\begin{displaymath}
a_{\nu \mu}\;:=\;\sum_{j=1}^{n}C_{j}(\alpha_{\nu} \alpha_{\mu}),
\end{displaymath}
and $<.,.>$  is the standard scalar product.

The matrix $A_{\alpha}$ is a {\bf symmetric positive definite} matrix.
We can thus write
\begin{displaymath}
A_{\alpha}\;=\;B_{\alpha}^{2},
\end{displaymath}
for some {\bf symmetric matrix} $B_{\alpha}$. Therefore,
($B^{*}_{\alpha} = B_{\alpha}$)
\begin{displaymath}
<A_{\alpha}x,x>=<B^{2}_{\alpha}x,x>=<B_{\alpha}x,B^{*}_{\alpha}x>=
\mid\mid B_{\alpha}x \mid\mid^{2_{n}}
\end{displaymath}
and
\begin{displaymath}
exp(-\pi (N(A)^{2}d(k))^{-1/n}<A_{\alpha}x,x>)=exp(-\pi
(N(A)^{2}d(k))^{-1/n}\mid\mid B_{\alpha}x
\mid\mid^{2}_{n})=\omega_{B_{\alpha}}(x)\;;x \in \LR^{n}.
\end{displaymath}
Thus, $C(A)$ {\bf corresponds} to $B_{\alpha}$ in Lang's
considerations of this gaussian fixed point.
From [L, III.1] it is immediately follows that the
inverse matrix of $A_{\alpha}$ is given by
\begin{equation}
<A_{\alpha}^{-1}x, x >\;=\;\sum_{j=1}^{n}
\mid\sum_{v=1}^{n}C_{j}(\alpha_{v}^{*})x_{v}\mid ^{2}.
\end{equation}
Furthermore, the absolute value of the discriminant is
\begin{displaymath}
\mid D_{k}(\alpha_{1}, ..., \alpha_{n})\mid \;=\;det(A_{\alpha}).
\end{displaymath}

One can establish the value of $H_{r}^{0}(E(k))$ exactly in the same
way as in [L, XIII.3]. More exactly, it is not difficult to calculate
that
\begin{equation}
H_{r}^{0}(E(k))\;=\;2^{r_{1}-1}R(k),
\end{equation}
where $R(k)$ is the {\bf regulator} of $k$ defined as follows : let
$u_{1}, ... , u_{r_{1}+r_{2}}$ be {\bf independent generators} for the
unit group $U(k)$ (modulo roots of unity) (the {\bf Dirichlet's theorem}).
The absolute value of the determinant
\begin{equation}
det[N_{v}log\mid C_{j}(u_{v})\mid]
\end{equation}
(here $N_{v}$ - as usual - denotes a local degree) is independent of the
choice of our generators $\{u_{j}\}$ and is called the {\bf regulator}
$R(k)$ of the field $k$. We note that this regulator, like all
determinants, can be interpreted as a volume of a parallelotope in
$(r_{1}+r_{2})$-space.

Finally, the zeta function $\zeta_{k}(s)$ has a simple pole at $s=1$
with a {\bf residue} equal to
\begin{displaymath}
\frac{2^{r_{1}}(2\pi)^{r_{2}}h(k)R(k)}{w(k)\mid d(k) \mid}
\end{displaymath}
and the non-zero constant $\lambda_{k}$ in the {\bf zero-polar
factor} (trivial zeta) from the Face theorem
\begin{equation}
\lambda_{k}\;=\;\frac{2^{r_{1}}h(k)R(k)}{w(k)}.
\end{equation}
\end{re}

\section{RH-fixed points of ${\cal F}_{n}$}

In this section we present constructions which lead to the
derivation of the main technical tool of this paper - the {\bf harmonic
notion} of an RH-fixed point of the n-dimensional real Fourier
transform. We present here a more abstract and brief version of the
technique which was originally developed in [$M_{A}$] for the proof of the
Riemann hypothesis.

Let $V$ be a {\bf real vector space} endowed with an {\bf idempotent
endomorphism} $F :V \longrightarrow V$, i.e. $F^{2} = I_{V}$, where
$I_{V}$ denotes the identity endomorphism of $V$.

Let us consider the {\bf purely algebraic} notion of the {\bf
quasi-fixed point of $F$} associated with a parameter $l \in \LC$ and
an element $v \in V$ :
\begin{equation}
Q_{l}(F)(v)\;=\;Q_{l}(v)\;:=\;v\;+\;lF(v).
\end{equation}

Let us observe that if $l=1$ then $Q_{1}(v)$ is a {\bf fixed point} of
$F$, i.e.
\begin{equation}
F(Q_{1}(v))=F(v)+F^{2}(v)=F(v)+v =Q_{1}(v),
\end{equation}
and if $l=-1$ then $Q_{-1}(v)$ is a (-)fixed point of $F$, i.e.
\begin{displaymath}
F(Q_{-1}(v))=F(v)\;-\;F^{2}(v)\;=\;-(v\;-\;F(v))=-Q_{-1}(v).
\end{displaymath}

We obtain the following result on the existence of quasi-fixed points
\begin{lem}({\bf Existence of quasi-fixed points}).

For each $v_{0} \in V$ and $l \ne \pm 1$ the {\bf formula}
\begin{equation}
v_{l}\;:=\;\frac{v_{0}}{1-l^{2}}\;-\;\frac{lF(v_{0})}{1-l^{2}}
\end{equation}
gives the solution of the following {\bf Abstract Fox Equation}( AFE
in short , cf. also [$M_{A}$])
\begin{equation}
(AFE_{V}^{l})\;\;v_{l}\;+\;lF(v_{l})\;=\;v_{0}.
\end{equation}
\end{lem}

Lemma 2 shows that making a simple algebraic calculus, we cannot obtain a
{\bf singular} solutions of $AFE_{V}$, since the formula (3.59) {\bf has no
sense} for $l=\pm 1$.

Moreover, we see that on the ground of {\bf classical logic} the
$\pm$ fixed point $Q_{\pm 1}(v_{0})$  cannot be the solution of $(AFE_{V}^{\pm})$
, $Q_{\pm}(v_{\pm})=v_{0}$ if $v_{0}$ is not a $\pm$fixed point of $F$.

Let us denote the real subspace of $V$ of all $\pm$fixed points $v_{0}$ of $F$ in
$V$ by $Fix_{\pm}(F)$, i.e. $F(v_{0}) = \pm v_{0}$. We thus see that
the condition
\begin{equation}
v_{0}\;\in\; Fix_{\pm}(F)
\end{equation}
is a {\bf necessary condition} for the {\bf existence} of solutions
$Q_{\pm1}(v_{0})$ of $(AFE_{V})$.

We construct $Q_{\pm 1}(v_{0})$ using the {\bf averaging procedure} for the
family $\{Q_{l}(v_{0}) : l \ne \pm 1\}$, originally constructed in
[$M_{A}$].

As in [$M_{A}$] , it will be very convenient to use the unique non-commutative
field of {\bf Hamilton quaternions} $\LH$ (the {\bf Einstein space-time space}).

We denote a {\bf Haar measure} of the additive group $(\LH, +)$, by
$H_{\LH}$, i.e. the {\bf standard Lebesgue measure} $d^{4}h$ of the
vector space $\LR^{4}$ (the Einstein space-time).

For each $M,N >0$ we consider the {\bf hamiltonian segments} (rings)
\begin{equation}
S(M,N)\;:=\;\{h \in \LH : M \le \mid h \mid_{\LH} \le N\},
\end{equation}
where in all the sequel $\mid \cdot \mid_{\LH} := \mid\mid \cdot
\mid\mid_{4}$ is the standard Euclidean norm on $\LR^{4}$.

Finally, we consider the {\bf invertion $I_{\LH}$} of $\LH$
\begin{equation}
I_{\LH}(l)\;:=\;l^{-1}\;,\;l\in \LH^{*}:=\LH-\{0\}.
\end{equation}
Let us observe that $I_{\LH}$ is only a set-automorphism (and not a
group automorphism of the multiplicative group $\LH^{*}$, since it is
not commutative).

Each {\bf automorphism} $\lambda$ of $(\LH, +)$ changes the Haar
measure $H_{\LH}$ into $cH_{\LH}$  with $c \in \LR_{+}^{*}$ (the {\bf
von Neumann-Weil theorem}). The number $c$ does not depend on the choice
of Haar measure. It is denoted by $\Delta_{\LH}(\lambda)$ and is
called the {\bf Haar module} of $\lambda$. It is defined by any of the
equivalent formulas given below (cf.[W,I])
\begin{displaymath}
(W_{m})\;\;H_{\LH}(\lambda(B))\;=\;\Delta_{\LH}(\lambda)H_{\LH}(B)
\end{displaymath}
or
\begin{displaymath}
(W_{i})\;\;\int f(\lambda^{-1}(x))dH_{\LH}(x)=\Delta_{\LH}(\lambda)\int
f(x)dH_{\LH},
\end{displaymath}
where $B$ is any Borel set and $f$ is any integrable function with
$\int f dH_{\LH} \ne 0$.

The second formula can be symbolically written in the form:
\begin{displaymath}
dH_{\LH}(\lambda(x))\;=\;\Delta_{\LH}(\lambda)dH_{\LH}(x).
\end{displaymath}

If $h \in \LH^{*}$ is arbitrary, then the formula : $M_{h}(x) := h\cdot
x, x \in \LH$ defines a linear multiplication automorphism of $(\LH,
+)$. We set
\begin{displaymath}
\Delta_{\LH}(h)\;:=\; \Delta_{\LH}(M_{h})\;,\;h \in \LH^{*},
\end{displaymath}
and moreover, we define $\Delta_{\LH}(0) := 0$. It is well-known that
(cf.e.g. [W, I.2 and Corrolary 2])
\begin{equation}
\Delta_{\LH}(h)\;=\;\mid h \mid^{4}_{\LH}\;=\;\mid\mid h
\mid\mid_{4}^{4}.
\end{equation}

We denote the {\bf invertion} of $\LH^{*}$ by $I_{\LH}(h) := h^{-1}, h
\in \LH^{*}$. Unfortunately, $I_{\LH} =I$ is not a group automorphism of
$\LH^{*}$ , since $\LH^{*}$ is not commutative! However, it is still a very
{\bf crucial topologically-algebraic} map of $\LH^{*}$ of order $2$ :
$I_{\LH}^{2} = id_{\LH}$.

Thus, beside such an important invariant of $\LH$ like the Galois
group $Gal(\LH/\LR)$, we have an additional important {\bf invariant} of
$\LH$ - the {\bf invertion group}  $Inv(\LH^{*}) := \{id_{\LH^{*}},
I_{\LH^{*}}\}$ of $\LH^{*}$ (cf. [$M_{A}$]).

It is well-known (cf.e.g.[$M_{A}$, Lem.4])  that
\begin{equation}
dH_{\LH^{*}}(h)\;:=\;\frac{dH_{\LH}(h)}{\mid h \mid_{\LH}^{4}}\;,\;h
\in \LH^{*}.
\end{equation}
is a (left) {\bf Haar measure} of the multiplicative group $\LH^{*}$.
Moreover, it would be convenient to recall the {\bf simple algebraic
-measure formulas} for $H_{\LH}$ and $H_{\LH^{*}}$ given below (cf.
$M_{A}$, Prop.3) : for each {\bf integrable} function $f$ on $\LH^{*}$
we have :
\begin{equation}
\int_{\LH^{*}}f(h^{-1})dH_{\LH^{*}}(h)
\;=\;\int_{\LH^{*}}f(h)dH_{\LH^{*}}(h) .
\end{equation}
i.e. $H_{\LH^{*}}$ is the {\bf invariant measure} (or the {\bf
Bogoluboff-Kriloff measure}) of the {\bf dynamical system} $(\LH^{*},I_{\LH})$.
Moreover,
\begin{equation}
\int_{\LH^{*}}f(h)dH_{\LH}
\;=\;\int_{\LH^{*}}\frac{f(h^{-1})dH_{\LH}(h)}{\mid h \mid_{\LH}^{8}}.
\end{equation}

For each $N > M >0$ we consider the {\bf compact $\LH$-rings}
\begin{equation}
R_{\LH}(M,N)\;:=\;\{h \in \LH : M \le \mid h \mid_{\LH} \le N\},
\end{equation}
and the corresponding {\bf dynamical sub-system} of $(\LH^{*}, I_{\LH})$
\begin{equation}
D_{\LH}(M,N)\;:=\;(R_{\LH}(M,N) \cup R_{\LH}(N^{-1},M^{-1}), I_{\LH}),
\end{equation}
with $M, N>1$.

From (3.66) we immediately obtain that the formula
\begin{equation}
\beta_{\LH}(A)\;:=\;\int_{A}\frac{d^{4}h}{\mid 1- h^{2}\mid ^{4}}\;;\;h
\in D_{\LH}(M,N) ,
\end{equation}
gives an {\bf invariant measure} of $D_{\LH}(M,N)$ - the {\bf Herbrand
distribution} of $\LH^{*}$ (cf.[$M_{A}$]). In particular, the measure $\beta
_{\LH}$ satisfies the condition
\begin{equation}
\beta_{\LH}(I_{\LH}^{-1}(A))\;=\;\beta_{\LH}(A).
\end{equation}

We use below the theory of the {\bf sextet} $(\LH, \Delta_{\LH}, H_{\LH},
R_{\LH}(M,N), \beta_{\LH}, R_{\LH})$ .

For the sake of completness, we also briefly recall here two deep
and difficult results from {\bf analytic potential theory} of $\LH$
explored in [$M_{A}$] :

(1) {\bf The Riesz theorem} (cf.[HK, Sect.3.5, Th.3.9]).

Let $s = s(x)$ be a {\bf subharmonic} function in a domain of
$\LR^{6}$. Then there exists a {\bf hamiltonian Riesz measure} $R_{\LH}$
and a {\bf harmonic function} $h(x)$ outside a compact set $E$, such that
\begin{displaymath}
(RT)\;s(x)\;=\;\int_{E}\frac{dR_{\LH}(y)}{\mid\mid x-y
\mid\mid_{6}^{4}}\;+\;h(x)\;;\;x \in \LR^{6}.
\end{displaymath}
(2) {\bf Brelot's theorem} (cf.[HK, Sec.36, Th3.10] - on the
existence of {\bf harmonic measures}).

Let $D$ be a {\bf regular} and {\bf bounded domain} of $\LR^{n}$ with
border $\partial D$. Then, for each $x \in D$ and arbitrary Borel
set $B$ of $\partial D$ , there exists a {\bf unique number}
$\omega(x,B:D)$, which is a {\bf harmonic function in $x$} and {\bf
probabilty measure in $B$} and moreover, for each {\bf semicontinuous
function} $f(\xi)$ on $\partial D$ the formula
\begin{equation}
(DP)\;\tilde{f}(x)\;=\;\int_{\partial D}f(\xi)d\omega(x,\xi;D)\;;\;x \in
D-\partial D ,
\end{equation}
gives the {\bf harmonic extension} of $f$ from $\partial D$ to $D$.

The family of harmonic measures $\omega(D) :=\{\omega(x,\cdot;D): x \in
D\}$  solves the {\bf Dirichlet problem}(DP) for a pair $(D, \partial D)$
and if a solution exists it is {\bf unique}.

In [$M_{A}$] we introduced the following formal definition of the {\bf
Abstract Hodge Decomposition}: let $f : X \longrightarrow \LC$ be a
function and $K : X \times I \longrightarrow \LC$ another "kernel"
function. A measure $H_{f}$ on a $\sigma$-field of subsets of $I$  gives
the Abstract Hodge Decomposition of $f$, if the following integral representation
is satisfied
\begin{displaymath}
(AHD_{f}) \;f(x)\;=\;\int_{I}K(x,i)dH_{f}(i) \;;\;x \in X .
\end{displaymath}
We call the measure $H_{f}$, which appeares in $(AHD_{f})$ the {\bf Hodge
measure} of $f$.
\begin{pr}({\bf Existence of $AHD_{\LH}$}).

There exists such a Borel probability measure $R_{\LH}$ (a {\bf
hamiltonian Riesz measure}) on the 3-dimensional sphere $S^{3} := \{h
\in \LH : \mid h \mid_{\LH}=1\}$, such that for each
\begin{displaymath}
r \in X_{\LH}(M,N)\;:=\;R_{\LH}(M,N)\cup R_{\LH}(N^{-1},M^{-1})
\end{displaymath}
with $N>M>1$, the following abstract Hodge decomposition ($AHD_{\LH}$
in short) holds :
\begin{displaymath}
(AHD_{\LH})\;\Delta^{-1}_{\LH}(r^{2})\;=\;\int_{S^{3}}\frac{dR_{\LH}(h)}
{\Delta_{\LH}(r^{2}-h^{2})}\;=\;{\cal H}_{\LH}(R_{\LH})(r).
\end{displaymath}
\end{pr}
{\bf Proof}. Let $\epsilon_{n} >0$ be an arbitrary sequence, which
converges to zero. Then the functions $\mid\mid \cdot
\mid\mid_{6}^{-(4+\epsilon_{n})}$  are {\bf subharmonic} ( as suitable powers
of a {\bf harmonic function}) and obviously they are {\bf not
harmonic}! Therefore, according the {\bf Riesz theorem}, there exists a
sequence of {\bf Riesz measures} $\{R_{n}\}$ and a sequence $\{h_{n}\}$
of harmonic functions {\bf inside} of $S^{3}$ with the property
\begin{equation}
\mid\mid r \mid\mid_{6}^{-(4+\epsilon_{n})}
\;=\;\int_{S^{3}}\frac{dR_{n}(h)}{\mid\mid r-h
\mid\mid_{6}^{4}}\;+\;h_{n}(r).
\end{equation}
Since $dR_{n}(x) =\nabla(\mid\mid x \mid\mid_{6}^{(4+\epsilon_{n})})dx$
(cf. [HK, Section 3.5]) , the sequence $\{R_{n}(S^{3})\}$ is
{\bf bounded}, i.e. $R_{n}(S^{3}) \le A$, for some $A>0$ and all $n \in
\LN$.

According to {\bf Frostman's theorem} (cf. [HK, Theorem 5.3]), we
can choose a subsequence $\{R_{n_{p}}\}$, which is {\bf weakly
convergent} to a limit measure $R_{\infty}$ on $S^{3}$, i.e. $R_{\infty}
:= (w)\lim_{p\longrightarrow \infty}R_{n_{p}}$ and
\begin{equation}
\mid\mid r
\mid\mid_{6}^{-4}\;=\;\int_{S^{3}}\frac{dR_{\infty}(x)}{\mid\mid r-x
\mid\mid_{6}^{4}}\;+\;h(r).
\end{equation}

On the other hand, according to {\bf Brelot's theorem} applied to
the triplet $(B^{6}(1),0,S^{5})$ - there exists a {\bf harmonic
measure} $\omega(\cdot):=\omega(\cdot,0,B_{6})$ on $S^{5}$ with the
property

\begin{equation}
\mid\mid r \mid\mid_{6}^{-4}\;=\;\int_{S^{5}}\frac{d\omega(y)}{\mid\mid
r-y\mid\mid_{6}^{4}}\;,\;r \in (B^{6})^{c}.
\end{equation}

Let us denote the natural inclusion by $j_{35} : S^{3} \longrightarrow
S^{5}$; $j_{35}(h)\;=\; (h,0,0)$. Then (3.74) can be written of the
form
\begin{equation}
\mid\mid r \mid\mid_{6}^{-4}\;=\;\int_{S^{5}}
\frac{d(j_{35}^{*}R_{\infty})(y)}{\mid\mid y -r
\mid\mid_{6}^{4}}\;+\;h(r).
\end{equation}
Let us consider the {\bf continuous function} $f(\xi) := \mid\mid \xi
\mid\mid_{6}^{-4}$ on $S^{5}$ and (for a while) take $\mu$ to be one of the
two measures : $j_{35}^{*}(R_{\infty})$ or $\omega$. Finally, let us
consider the potential $\int_{S^{5}}\frac{d\mu(y)}{\mid\mid
r-y\mid\mid_{6}^{4}}$. The vectors $r, y \in \LR^{6}$  can be
considered as {\bf Cayley numbers} from $\LR^{8} = \LH \times \LH$
and $\mid\mid \cdot \mid\mid_{6}$  as the restriction of the {\bf
Cayley norm}. Since Cayley numbers form a non-commutative and
non-associative algebra with {\bf division}, then we can write
\begin{displaymath}
\int_{S^{5}}\frac{d\mu(y)}{\mid\mid
r-y\mid\mid_{6}^{4}}=\int_{S^{5}}\frac{d\mu(y)}{\mid\mid y(1-r/y) \mid\mid_{6}^{4}}
=:\int_{S^{5}}\frac{d\nu(r,y)}{\mid\mid y
\mid\mid_{6}^{4}}=\int_{S^{5}}f(\xi)d\nu(r,\xi).
\end{displaymath}
Thus, both the formulas (3.75) and (3.76) give the solution of the
{\bf Dirichlet problem} for $((B^{6})^{c},S^{6},f)$. From the {\bf
uniqueness} of the solution of the Dirichlet problem (cf. [HK,
Th.1.13]), we obtain:
\begin{displaymath}
j_{35}^{*}(R_{\infty})\;=\;\omega \;\;and\;\; h\equiv 0\;on\;
(B^{6})^{c},
\end{displaymath}
since $h$, as a difference between a harmonic function and a potential, is
also harmonic on $\LR^{6}- S^{5}$.
Hence, restricting ourselves in (3.76) for $r = h \in \LH$, we finally
obtain
\begin{equation}
\Delta_{\LH}(h^{-1})\;=\;\int_{S^{3}}\frac{dR_{\infty}(x)}{\Delta_{\LH}(h
-x)} \;,\;h \in S_{\LH}(M,N).
\end{equation}

Let us consider a {\bf branch of the hamiltonian square root} $\sqrt{\cdot}$
and the induced map of measure spaces : $\sqrt{\cdot} :
(S^{3},R_{\infty})\longrightarrow (S^{3}, \sqrt{\cdot}^{*}R_{\infty})$,
substituting $h^{2}$ for $h$ and $R_{\LH}$  for
$\sqrt{\cdot}^{*}R_{\infty}$ we obtain the above proposition.

We will use the {\bf hamiltonian sextet} (from analytic potential theory)
\begin{displaymath}
(\LH, \mid \cdot \mid_{\LH}, \Delta_{\LH}, H_{\LH},
H_{\LH^{*}},R_{\LH}, \beta_{\LH})
\end{displaymath}
in the averaging procedure given below to obtain, singular solutions of
$(AFE_{V})$.

A similar result is much easier to obtain using the completely
different nature of locally compact rings - the small adeles , i.e.
working with the {\bf adic potential theory}.

As we will show below, in the p-adic case, the required {\bf algebraic}
potential theory is simpler, in opposite to the strongly analytic
potential theory of $\LR^{m}$. Therefore, the p-adic fields (and
generally local non-archimedean fields are - in such a way we see them
today - are missing links - to the needed maths constructions).

Let $H_{p}$ denotes the {\bf Haar measure} of the additive group of the
p-adic number field $\LQ_{p}$. The main reason that the algebraic
potential theory over $\LQ_{p}$ is simpler that the analytic one over
$\LR^{m}$ is the quite different behaviour of Haar measures on totally
- disconnected fields with compare to Haar measures on the connected
fields. In particular, $\LZ_{p}^{*}$ is {\bf open}, and therefore
$H_{p}(\LZ_{p}^{*}) \ne 0$, whereas in the case of {\bf connected}
local fields $K$ we have
\begin{displaymath}
H_{K}(S_{K})\;=\;0,
\end{displaymath}
where $S_{K}$ is the unit sphere in $K$.

Thus, it is convenient to normalizeed $H_{p}$ in such a manner that
\begin{displaymath}
H_{p}(\LZ_{p}^{*}=S_{p})\;=\;(1\;-\;p^{-1}),
\end{displaymath}
(the Euler component in $\zeta_{\LQ}^{-1}(1)$)

Let $p$ and $q$ be two different {\bf prime numbers} and $\LQ_{p}$ and
$\LQ_{q}$ be the fields of p-adic and q-adic numbers. For the
convenience we take the non-canonical choice of $\mid . \mid_{q}$  by
putting $\mid q \mid_{q}= 1/p$!

We denote by $\LQ_{pq}$ the product $\LQ_{p}\times \LQ_{q}$, being the
locally compact abelian ring with a large set of invertible elements (its
completion has Haar measure zero). We denote its {\bf Haar module} by
$\Delta_{pq}$, and its {\bf Haar measure} by $H_{pq}$. $R_{pq}$ is the
{\bf adic Riesz measure}, $I_{pq}$ is the {\bf invertion automorphism} of
$\LQ_{pq}^{*}$ and finally, we denote the {\bf adic Herbrandt distribution}
(the Bogoluboff-Krilov measure) of $\LQ_{pq}$ by $\beta_{pq}$. Then
the sextet
\begin{displaymath}
(\LQ_{pq},\LQ_{pq}^{*},\Delta_{pq},H_{pq},I_{pq},\beta_{pq})
\end{displaymath}
enables us to show, in a relatively simple, algebraic way,  the existence
of the below {\bf $pq$-adic Abstract Hodge Decomposition}($(AHD_{pq})$
in short).

According to the Weil's Lemma 2 (see [W, I.2., Lemma 2]) we have
\begin{displaymath}
\Delta_{pq}(x)\;=\;mod_{\LQ_{pq}}(x)\;=\;\mid x_{p} \mid_{p} \mid x_{q}
\mid_{q}\;,\;x=(x_{p},x_{q}) \in \LQ_{pq}^{*}.
\end{displaymath}

We define a {\bf sub-dynamical system} $D_{pq}(M,N) = (X_{pq}(M,N),
I_{pq})$ of the dynamical system of the small adeles $(\LQ_{pq},
I_{pq})$. Moreover, in the sequel we simply write $D_{pq}$ and $X_{pq}$
instead of $D_{pq}(M,N)$ and $X_{pq}(M,N)$, respectively. The compact
topological space $X_{pq}$ is defined as follows: let $M, N \in \LN^{*}$
be such that $1 \le M < N$. Then
\begin{displaymath}
X_{pq}(M,N)\;:=\;\{x \in \LQ_{p} \times \LQ_{q}: x=(x_{p},x_{q}), \mid
x_{q} \mid_{q}=1, \mid x_{p} \mid_{p}\in [p^{-N+1},p^{-M}]\cap
[p^{M},p^{N-1}]=
=\;I_{\LR}([p^{M},p^{N-1}])\cap [p^{M},p^{N-1}]\}.
\end{displaymath}
Finally, let us consider the the p-dic {projection}
$P_{p}:\LQ_{pq}\longrightarrow \LQ_{p}, P_{p}(x_{p},x_{q})=x_{p}$ and
$I_{p}$-{\bf invariant} function
\begin{displaymath}
{\cal I}_{p}(\lambda)\;:=\;\frac{\mid \lambda \mid_{p}}{\mid
1\;-\;\lambda^{2}\mid_{p}}\;,\;\lambda \in \LQ_{p}^{*}-\{1\}.
\end{displaymath}
Under the above notations we have
\begin{lem}({\bf On the pq-adic Herbrandt measure $d\beta_{pq}$}).

The formula
\begin{displaymath}
d\beta_{pq}(x_{p},x_{q}):=\frac{d(H_{p}\times
H_{q})(x_{p},x_{q})}{\mid 1\;-\;x_{p}^{2}\mid_{p}} = \frac{{\cal
I}_{p}(P_{p}(x_{p},x_{q}))d(H_{p}\times
H_{q})(x_{p},x_{q})}{\Delta_{pq}((x_{p},x_{q}))}
\end{displaymath}
gives a {\bf Bogoluboff-Kriloff measure}( {\bf Herbrandt distribution})
of $D_{pq}$.
\end{lem}
{\bf Proof}. Since $\Delta_{pq}$ is the Haar module of $\LQ_{p}$ (like
$\Delta_{k_{A}^{*}}(z)=\prod_{v \in P}\mid z_{v} \mid_{v}$) in the case
of {\bf ideles} (see e.g. [W] and [Ko]), then the equality
\begin{displaymath}
\Delta_{pq}(x_{p},x_{q})\;=\;\mid x_{p} \mid_{p},
\end{displaymath}
holds on $X_{pq}(M,N)$. Hence we get
\begin{displaymath}
\frac{d(H_{p}\times H_{q})(x_{p},x_{q})}{\mid 1-
x_{p}^{2}\mid_{p}}=\frac{\mid x_{p}\mid_{p}}{\mid 1- x_{p}^{2}
\mid_{p}}\cdot \frac{d(H_{p}\times H_{q})(x_{p},x_{q})}{\mid
x_{p}\mid_{p}}=
\end{displaymath}
\begin{displaymath}
{\cal I}_{p}(P_{p}(x_{p},x_{q}))\cdot \frac{d(H_{p}\times H_{q})(x_{p}
, x_{q})}{\Delta_{pq}(x_{p},x_{q})}.
\end{displaymath}
Since $\frac{d(H_{p}\times H_{q})}{\Delta_{pq}}$ is a
Bogoluboff-Kriloff measure of $(\LQ_{pq}^{*}, I_{pq}), P_{p} \circ
I_{pq}=I_{p}\circ P_{p}$ and ${\cal I}_{p}$ is $I_{p}$-invariant, that
we really see that the above formula gives a $pq$-adic
Bogoluboff-Kriloff measure of $D_{pq}$. (Let us remark the importance
of the fact that $H_{q}(\LZ_{q}^{*}) \ne 0$).

\begin{pr}({\bf The existence of $AHD_{pq}$}).

\begin{displaymath}
(AHD_{pq})\;\;\Delta_{pq}(x^{-2})\;=\;\int_{S_{pq}}\frac{dR_{pq}(y)}
{\Delta_{pq}(x^{2}-y^{2})} \;\;if\;x \in X_{pq}(M,N) ,
\end{displaymath}
where $S_{pq}:=\{x \in \LQ_{pq} : \Delta_{pq} =1\}$ is the unit adic
sphere.
\end{pr}
{\bf Proof}. Let $x \in X_{pq}(M,N)$ be arbitrary. Then $x = (x_{p},x_{q})$
with $\mid x_{q}\mid_{q}=1$ and therefore $\Delta_{pq}(x) = \mid x_{p}
\mid_{p}$ ( we can identity the p-adic field $\LQ_{p}$ with the subset
$\LQ_{p}\times \{1\}$ of $\LQ_{pq}$). But $\mid x_{p} \mid_{p} \ge
p^{-N+1}$ , and therefore according to the {\bf ultrametricity} of
$\mid \cdot \mid_{p}$ (see e.g. [W, I.2., Corrolary 4]) for all $ y \in
\LQ_{p}$ with $\mid y \mid_{p} = p^{-N}$ we have
\begin{displaymath}
\mid x_{p}^{2} \mid_{p}\;=\; \mid x_{p}^{2}\;-\;y^{2} \mid_{p}.
\end{displaymath}
Integrating the both sides of the inverse of the above equality with
respect to the Haar measure $H_{p}$ on $p^{N}\LZ_{p}^{*}$, for each
$\eta \in \LQ$ with $\mid \eta \mid_{p}=1$, we obtain :
\begin{displaymath}
\frac{1}{\mid x_{p}^{2}
\mid_{p}}\;=\;\frac{1}{H_{p}(p^{N}\LZ_{p}^{*})}\int_{p^{N}\LZ_{p}^{*}}
\frac{dH_{p}(\xi)}{\mid x_{p}^{2}\;-\;(\eta \xi)^{2} \mid_{p}}.
\end{displaymath}
Let us denote : $d \nu_{p}(\xi) :=
\frac{dH_{p}(\xi)}{H_{p}(p^{N}\LZ_{p}^{*})}$. Then, for all $\eta \in
\LQ$ with $\mid \eta \mid_{p} = 1$, the above equality can be written as
\begin{displaymath}
\frac{1}{\mid x_{p}^{2}\mid_{p}}\;=\;\int_{p^{N}\LZ_{p}^{*}}\frac{d
\nu_{p}(\xi)}{\mid x_{p}^{2}\;-\;(\eta \xi)^{2} \mid_{p}}.
\end{displaymath}
(we non-standartly assumed that $\mid q \mid_{q} = p^{-1})$.

Let $F$ be any {\bf finite} subset of $\{\eta \in \LQ : \mid \eta
\mid_{p}=1, \mid \eta \mid_{q} = p^{-N}\} \subset p^{-N}\LZ_{q}^{*}$.
For an arbitrary subset $A$ of $\LQ$ we define the measure $\mu_{q}^{F}$
by
\begin{displaymath}
\mu_{q}^{F}(A)\;:=\;\sum_{f \in F} \frac{\delta_{f}(A)}{\#F} ,
\end{displaymath}
where $\delta_{f}$ is the Dirac measure at $f$.

Let us consider the measure $(\nu_{p} \times \mu_{q}^{F})$. Summing the
both sides of the previous measure representation on $F$ we obtain
\begin{displaymath}
\frac{1}{\mid x_{p}^{2} \mid_{p}}\;=\;\int \int_{p^{N}\LZ_{p}^{*}\times
p^{-N}\LZ_{q}^{*}\subset S_{pq}}\frac{d(\nu_{p}\times
\mu_{q}^{F})(\xi,\eta)}{\mid x_{p}^{2}\;-\;(\eta \xi)^{2}\mid_{p}}.
\end{displaymath}
Let us look at the natural inclussion $j_{pq} : p^{N}\LZ_{p}^{*}\times
p^{-N}\LZ_{q}^{*}$ as on a {\bf random variable}
$j_{pq}:(p^{N}\LZ_{p}^{*}\times p^{-N}\LZ_{q}^{*},\nu_{p}\times
\mu_{q}^{F})\longrightarrow S_{pq}$. Then the {\bf distribution} $R_{pq}$
of $j_{pq}$ we will be called the {\bf (p,q)-adic Riesz measure} and
the right-hand side of the above formula we can finally write in the
form :
\begin{displaymath}
\frac{1}{\mid x_{p}^{2} \mid_{p}}\;=\;\int_{S_{pq}}\frac{R_{pq}(dy)}{\mid
x_{p}^{2}\;-\;P_{p}(y)^{2} \mid_{p}}.
\end{displaymath}
Combining the above formulas, we obtain the proof of the existence of
the $(AHD_{pq})$. It also shows that the proof of $(AHD_{pq})$ is
possible in a completely {\bf algebraic way}.
\begin{re}
Probably the first mathematician, who considered and applied the p-adic
potential theory was {\bf Kochubei}. In the case of p-adic fields
$\LQ_{p}$ , the $\LQ_{p}$-Hilbert transforms probably first were
considered  in the {\bf Vladimirov}  et al.'s paper [VWZ] as the
$\gamma$-order {\bf derivative} $D^{\gamma}f$ of a locally constant
function $f$. It is describable by pseudo-differential operator and
explicitly written as
\begin{displaymath}
D^{\gamma}f(x)=\int_{\LQ_{p}}\mid \xi
\mid^{\gamma}_{p}\hat{f}(\xi)\chi_{p}(-\xi x)
H_{p}(d\xi)=\frac{p^{\gamma-1}}{1-p^{-\gamma-1}}\int_{\LQ_{p}}\frac{f(x)-
f(y)}{\mid x-y\mid_{p}^{\gamma+1}}H_{p}(dy),
\end{displaymath}
where $\chi_{p}$ is the additive character of $\LQ_{p}$ and
$\hat{f}(\xi)$ stand for the Fourier transformation
$\int_{\LQ_{p}}\chi_{p}(\xi x)f(x)H_{p}(dx)$ of a function $f$. A
deeper analysis of p-adic fractional differentiation $D^{\gamma}$ is
given in the Kochubei's book [Ka], where using {\bf Minlos-M\c{a}drecki's
theorem}, he established the existence of a {\bf Kochubei-Gauss measure}
$\mu$ over infinite-dimensional field extensions $\Omega_{p}$ of
$\LQ_{p}$, which is a {\bf harmonic measure} for $D^{\gamma}$ and
solves {\bf p-adic integral equations of a profile of wing of a plane}
in the case of $\Omega_{p}$ (see [Ka, Prop.6]).

The importance of $R_{pq}$ is also underline by the fact that
unfortunately, firstly we have the following negative result concerning
$(AHD_{pq})$.

{\bf Non-existence of solutions of p-adic profile of a wing in
functions}

Let $p$ be an arbitrary prime number. There is not exist an {\bf
absolutely continuous measure} $h_{p}$ w.r.t. the Haar measure $H_{p}$ (
the p-adic harmonic measure), which gives the following p-adic Abstract
Hodge Decomposition (cf. [$M_{A}$]) with the property
\begin{displaymath}
(AHD_{p})\;\frac{1}{\mid x
\mid_{p}^{2}}\;=\;\int_{\LZ_{p}^{*}}\frac{dh_{p}(y)}{\mid x^{2}-y^{2}
\mid_{p}}\;,\;x \in S_{p}(M,N).
\end{displaymath}
{\bf Proof}. The proof is based on the remarkable property of the Haar
measure $H_{p}: H_{p}(\LZ_{p}^{*})=1-p^{-1}$ (the {\bf Euler factor in
the Riemann zeta}). Assume (a contrary), that there exists a measure
$h_{p}$, which is absolutely continuous w.r.t. $H_{p} : h_{p}\ll
H_{p}$. Let us denote its density by $\omega_{p}$, i.e.
\begin{displaymath}
\omega_{p}\;=\;\frac{dh_{p}}{dH_{p}}.
\end{displaymath}
This conjecture permits us to apply the big and well-known machinery of
p-adic Fourier analysis to the problem of the existence of $(AHD_{p})$.
Reely, the $(AHD_{p})$ is obviously equivalent to the formula
\begin{displaymath}
\chi_{S_{p}(M,N)}(x)\mid x \mid_{p}^{-1} =\int_{\LQ_{p}}\mid x-y
\mid_{p}^{-1}\chi_{\LZ_{p}^{*}}(y)\omega_{p}(y)dH_{p}(y)\;=\;
\end{displaymath}
\begin{displaymath}
=\int_{\LQ_{p}}(\chi_{\LZ_{p}^{*}}\cdot \omega_{p})(x-z)\mid z
\mid_{p}^{-1}dH_{p}(z):=[(\chi_{\LZ_{p}^{*}}*\mid \cdot
\mid_{p}^{-1})](x)\;,\;x \in \LQ_{p},
\end{displaymath}
where $\chi_{A}$ denotes the characteristic function of a set $A$ and
$*$ means the p-adic {\bf convolution}.

If we apply the p-adic Fourier transform ${\cal F}_{p}$ to the both
sides of the above equalities, then we obtain
\begin{displaymath}
\hat{(\chi_{S_{p}(M,N)}\cdot \mid \cdot
\mid_{p}^{-1})}(\xi)\;=\;\hat{\mid \cdot \mid_{p}^{-1}}(\xi)\cdot
\hat{\chi_{\LZ_{p}^{*}}\cdot \omega_{p}}(\xi)\;,\;\xi \in \LQ_{p}.
\end{displaymath}
Let us observe that
\begin{displaymath}
\hat{\mid \cdot \mid_{p}^{-1}}(\xi)\;=\;\frac{(1-p)log \mid \xi \mid_{p}
}{plogp}\;,\;\xi \in \LQ_{p},
\end{displaymath}
see [Ka, Sect. 1.5, formula (1.29)]. Thus, the last equality is not
possible.

In the light of the above presented negative result, the previous above
result - on the existence of $(AHD_{pq})$ - gathers a greater value.

The small pq-adele ring $\LQ_{pq}$ is only one representant from a whole class
of "models", of the very similar nature, which can be used in the same
context.

(1){\bf The pq-adic vector space $\LQ_{[pq]}$}.

It is well-known (see e.g. [La, Sect.1]) that the "world" of {\bf
valuations} (or {\bf absolute values} or {\bf points}) is very rich. In
particular, we saw, how effective was the action of the defined below
pre-valuations $v_{pq}$, which gives $(AHD_{pq})$ of the p-adic
valuation $\mid \cdot \mid_{p}$ in the simple {\bf algebraic way}, if
we compare it with a difficult {\bf analytic} proof of $(AHD_{\LH})$ of
$\Delta_{\LH}^{-2}$. With a similar situation we have deal in the
famous Faltings' proof of the Mordell-Shafarevich-Tate conjectures. He
used so called {\bf heights of global fields}, which are some functions
defined by p-adic valuations (cf.e.g. [La, Fa]).

According to the principal theorem of the arithmetics each non-zero
rational number $x$ we can uniquely write of the form:
\begin{displaymath}
x\;=\;\frac{a}{b}p^{m}q^{n},
\end{displaymath}
where $a, b, m, n \in \LZ, (a,b)=1$, and $pq$ do not divide $ab$. For
such a rational $x$ we put
\begin{displaymath}
\alpha_{p}(x)\;:=\;m,\;\;\alpha_{q}(x)\;:=\;n,
\end{displaymath}
and
\begin{displaymath}
v_{pq}(x)\;:=\;p^{-(m+n)}\;=\;p^{-(\alpha_{p}(x)+\alpha_{q}(x))}.
\end{displaymath}
The functions $\alpha_{p} : \LZ \longrightarrow \LZ$ defined below are
called {\bf exponents} corresponding to $p$ and satisfies few simple
and nice elementary properties. A. Ostrowski showed their most
surprising properties : they are {\bf unique} arithmetical functions
(up to a constant - like Haar measures), which satisfy the five
mentioned above their elementary properties (see e.g. [Na1,Th.1.7(i)-(v)]).
In particular, the functions $v_{pq}$ satisfies the following condition
:
\begin{displaymath}
v_{pq}(x)\;=\;\mid x \mid_{p} \mid x \mid_{q},
\end{displaymath}
i.e. $v_{pq}$ has only one good "residual" above {\bf multiplicative}
property. Moreover
\begin{displaymath}
v_{pq}(x+y)\le max\{ v_{pq}(x), v_{pq}(y), \mid x \mid_{p}\mid y
\mid_{p},\mid y \mid_{p}\mid x \mid_{q}\},
\end{displaymath}
where $\mid \cdot \mid_{p}$ and $\mid \cdot \mid_{q}$ are p-adic and
q-adic valuations, respectively, but with the additional assumption
that
\begin{displaymath}
\mid q \mid_{q}\;=\;p^{-1}.
\end{displaymath}
Thus, $v_{pq}$ has a bad linear (ring) algebraic properties. In
particular, $v_{pq}(\cdot)$ is not a {\bf valuation} but only - let us
say - a {\bf pre-valuation}. Therefore, the {\bf completion} of $\LQ$
w.r.t. the metric type function : $d_{pq}(x,y)\;=\;v_{pq}(x-y)$ is
rather a pathological object and in particular, it is not a topological
field.

(2) By $\LQ_{(pq)}$ we donote the set $\{0,1, ... , pq-1\}(p,q)$ of all
{\bf double formal Laurent series} with coefficients in the set $\{0,1,
... ,pq-1\}$. Thus, each element $x$ of $\LQ_{(pq)}$ has the form :
\begin{displaymath}
x\;=\;\sum_{m=M}^{\infty}\sum_{n=N}^{\infty}a_{mn}p^{m}q^{n}\;;\;a_{mn}\
in \{0,1, ... , pq-1\}\;,\;M,N \in \LZ.
\end{displaymath}
If we establish a (non-canonical) ordering $<_{2}$ on the lattice
$\LZ^{2}$, in such a way that $(\LZ^{2}, \le_{2})$ and $(\LN, \le)$ are
isomorphic in the category of ordered sets : $h:(\LZ^{2}, \le_{2}) \simeq
(\LN, \le)$, and we establish the natural bijection
\begin{displaymath}
\LQ_{(pq)} \supset
x=\sum_{m=M}^{\infty}\sum_{n=N}^{\infty}a_{mn}p^{m}q^{n}\longrightarrow
\sum_{n=h(M,N)}a_{h(m,n)}p^{h(m,n)}\in \LQ_{p}\;=\;
\end{displaymath}
\begin{displaymath}
\;\;\;\;\;=\;\sum_{n=h(M,N)}a_{\pi(h(m,n))}p^{n},
\end{displaymath}
(where here $\pi$ denotes a respectible permutation of $\LN$), then we
can endow $\LQ_{(pq)}$ with the natural local field structure (transformed
from $\LQ_{p}$). Thus $\LQ_{(pq)}$ is a {\bf local field} isomorphic with
$\LQ_{p}$. Unfortunately, we cannot expect that $\LZ_{(pq)}$ is isomorphic
with $\LZ_{p}\otimes_{\LZ} \LZ_{q}$, where $\LZ_{(pq)}:=\{x \in
\LQ_{(pq)}: v_{pq}(x)\le 1\}$ (let us mention here that $\LF_{pq} \ne
\LF_{p}\otimes_{\LZ} \LF_{q}=0$).

To see the compactness of $\LZ_{(pq)}$ (and hence the local compactness
of $\LQ_{(pq)}$) it suffices to observe that the function $f$ from the
product $D$ of a countable many copies of the pq-elements set $\{0,1,
..., pq-1\}$ onto $\LZ_{(pq)}$ given by
\begin{displaymath}
f(\{a_{mn}\}_{m,n=0}^{\infty})\;=\;\sum_{m=0}^{\infty}\sum_{n=0}^{\infty}
a_{mn}p^{m}q^{n}
\end{displaymath}
is {\bf surjective} and {\bf continuous} in the Tichonov topology of
$D$. Since $D$ is compact then $\LZ_{(pq)}$ is compact as a continuous
image of the compact set.

On the other hand, we have natural inclusions (in the category of
sets): $i_{p} :\LQ_{p}\longrightarrow \LQ_{(pq)}$ and $i_{q}:
\LQ_{q}\longrightarrow \LQ_{(pq)}$.

Observe however, that with the multiplications defined as :
\begin{displaymath}
\alpha \cdot_{p} x\;:=\;i_{p}(\alpha)\cdot x \;\;and\;\;\alpha
\cdot_{q}x\;:=\;i_{q}(\alpha)\cdot x,
\end{displaymath}
(where $\cdot$ means the multiplication in $\LQ_{(pq)}$) $\LQ_{(pq)}$
{\bf is not a vector space} over $\LQ_{p}$ or over $\LQ_{q}$. Reely, if
it would be true, then obviously we would have:
$dim_{\LQ_{p}}\LQ_{(pq)}=+\infty$ and
$dim_{\LQ_{q}}\LQ_{(pq)}=+\infty$, what is impossible, since it is
well-known that LCA-vector spaces over local fields must be
finite-dimensional (see e.g. [W, I.2 Corrolary 2]). Moreover,
according to the {\bf Dantzing's description} of local fields, all
extensions of p-adic number fields $\LQ_{p}$, must be finite extensions
of such fields!
(3). {\bf The tensor product rings} $\LQ_{p}\otimes_{\LQ} \LQ_{q}$.

Let us observe that our main bi-adele (small pq-adele) ring
$\LQ_{pq}=\LQ_{p}\times \LQ_{q}$ is sufficiently good and "rich" for
our purpose, since from the point of view of the Haar-module theory the
set of its all {\bf non-invertible} elements :
$(\LQ_{pq}^{*})^{c}:=\LQ_{pq}-\LQ_{pq}^{*}$ is "small", i.e. its {\bf Haar
measure} is zero : $(H_{p}\times H_{q})((\LQ_{pq}^{*}))=0$.

For each $x \in \LQ_{pq}^{*}$ by $\Delta_{pq}(x)$ (or $mod_{\LQ_{pq}}(x)$)
we denoted the {\bf Haar module} of the {\bf automorphism} $x \longrightarrow
a \cdot x$ of $(\LQ_{p}\times \LQ_{q})^{+}$. Thus
\begin{displaymath}
\Delta_{pq}(x)\;=\;mod_{\LQ_{pq}}(x)\;:=\;\frac{(H_{p}\times
H_{q})(xX)}{(H_{p}\times H_{q})(X)},
\end{displaymath}
for arbitrary measurable set $X$ in $\LQ_{pq}$ with $0<(H_{p}\times
H_{q})(X)<+\infty$ (for example, for $X$ we can take any compact
neighbourhood of zero).

According to the Weil's Lemma 2 (see [W, I.2 , Lemma2]) we have
\begin{displaymath}
\Delta_{pq}(x)\;=\;mod_{\LQ_{pq}}(x)\;=\;\mid x_{p} \mid_{p}\mid x_{q}
\mid_{q}\;,\;x=(x_{p},x_{q})\in \LQ_{pq}^{*}.
\end{displaymath}

The above formula suggests that we can also descibe the pq-vectors from
$\LQ_{pq}$ in the terminology of the {\bf Grothendieck} tensor products.

Let us consider the {\bf algebraic} tensor product $\LQ_{p}\otimes_{\LQ}\LQ_{q}$
and the natural map $t_{pq}: \LQ_{p}\otimes_{\LQ}\LQ_{q}\longrightarrow
\LQ_{(pq)}$ defined by
\begin{displaymath}
t_{pq}(x \otimes y)\;=\;xy\;,\;x \in \LQ_{p},\;y \in \LQ_{q},
\end{displaymath}
(althought, according to the above mentioned troubles with the
multiplication in $\LQ_{(pq)}$ it is not {\bf algebraic}).

The Grothendieck $\pi$-norm $\mid \cdot \mid_{p}\otimes_{\pi}
\mid \cdot \mid_{q}$ (of $\mid \cdot \mid_{p}$ and $\mid \cdot
\mid_{q}$, see e.g. [MT]) is denoted in the sequel by $\mid \cdot
\mid_{pq}^{\pi}$ and is defined by
\begin{displaymath}
\mid x \mid_{pq}^{\pi}\;:=\;(\mid \cdot \mid_{p}\otimes \mid \cdot
\mid_{q})(x)\;:=\;
\end{displaymath}
\begin{displaymath}
\;=\;inf\{\sum_{i} \mid x_{i} \mid_{p}\mid y_{i}
\mid_{q}\;:\;x=\sum_{i}x_{i}\otimes y_{i}\;;\;x_{i} \in \LQ_{p},
y_{i}\in \LQ_{q}\}.
\end{displaymath}
The above $\pi$-norm is a {\bf cross-norm} (of $\mid \cdot \mid_{p}$
and $\mid \cdot \mid_{q}$), i.e.
\begin{displaymath}
(\mid \cdot \mid_{p}\otimes_{\pi} \mid \cdot \mid_{q})(x_{1}\otimes
x_{2})\;=\;\mid x_{1}\mid_{p} \cdot \mid x_{2} \mid_{q}\;;\;x_{1} \in
\LQ_{p}, x_{2}\in \LQ_{q}.
\end{displaymath}
Moreover
\begin{displaymath}
\mid x_{1} \otimes x_{2} \mid_{pq}^{\pi}=\mid x_{1}\mid_{p} \mid x_{2}
\mid_{q}= mod_{\LQ_{pq}}((x_{1},x_{2})) =
\end{displaymath}
\begin{displaymath}
=\;\Delta_{pq}((x_{1},x_{2}))
\end{displaymath}
$x_{1} \in \LQ_{p}, x_{2} \in \LQ_{q}$ and that norm is {\bf
archimedean}. By
$\LQ_{p}\hat{\otimes}_{\pi}\LQ_{q}=:\LQ_{pq}^{\otimes}$ we denote
the {\bf completion} of $(\LQ_{p}\otimes _{\LQ} \LQ_{q}, \mid \cdot
\mid_{pq}^{\pi})$.

Obviously $\LQ_{pq}^{\otimes}$ is a LCA-ring. Let us denote by $H_{pq}$
its {\bf Haar measure}.

Let $i_{pq} : \LQ_{p} \times \LQ_{q}\longrightarrow \LQ_{pq}^{\otimes}$
be the canonical inclusion homomorphism, i.e. $i_{pq}(x,y) = x \otimes
y$.

Since we have got {\bf Haar measures} $H_{p}$ and $H_{q}$ of
$\LQ_{p}^{+}$ and $\LQ_{q}^{+}$, respectively, then we can define their
{\bf tensor product} $H_{p}\otimes H_{q}$ on the LCA-subring $Im(i_{pq})$ :
\begin{displaymath}
i_{pq}^{*}(H_{p}\times H_{q})\;=:\;H_{p}\otimes_{\pi}H_{q}.
\end{displaymath}
It is easy to check that the tensor product $H_{p}\otimes_{\pi} H_{q}$
of Haar measures is a {\bf Haar measure} on $Im(i_{pq})$.

Thus, since the Haar measure $H_{pq}$ of $\LQ_{p}\hat{\otimes}_{\pi}\LQ_{q}$
is unique (up to a constant), then we can assume that
\begin{displaymath}
H_{pq}\;=\;H_{p}\otimes_{\pi}H_{q},
\end{displaymath}
i.e. the Haar measure $H_{pq}$ is the $\pi$-tensor product of the Haar measures
of $\LQ_{p}^{+}$ and $\LQ_{q}^{+}$ (see [MT]).

Finally, let us remark that the tensor products of gaussian measures in
Banach spaces were firstly considered by {\bf R. Carmona} and {\bf S. Chevet}
In [$M_{T}$] were defined and considered tensor products of
p-stable measures with $0<p<2$, in Banach spaces of stable type. The
above considered tensor products of Haar measures are also tensor
products of measures in the sense of the definition given in
[M$_{T}$].

(4){\bf The adele ring $\LQ_{\LA} \subset \prod_{p \in P}\LQ_{p}\times
\LQ_{\infty}(;=\LR)$ of $\LQ$}.

{\bf Ideles} were introduced by {\bf C. Chevalley} in [Ch] in 1936.
{\bf E. Artin} and {\bf G. Whaples} occured {\bf adeles} in [AW], where
they are called valuation vectors. The ring of adeles admits {\bf K.
Iwasawa's characterisation} (see [I]) in the following way : if $R$ is a
semi-simple commutative LCA-ring with a unit element, which is neither
compact nor discrete, and there is a field $K \subset R$, with the same
unit element, which is discrete and such that $R/K$ is compact, then $R$
is the {\bf ring of adeles} either over an algebraic number field or over an
algebraic function field with a finite fields of constants.

Topological properties of adeles and ideles were investigated by E.
Artin, K. Iwasawa, T. Tamagawa and J. Tate (see e.g.[N, Chapter VI]).

At the end of the ends, all the above considered versions of $\LQ_{pq}$
are closely related to each other and moreover we have the following
inclusion :
\begin{displaymath}
\LQ_{pq} \subset \LQ_{p}\otimes_{\LQ} \LQ_{q} \subset \LQ_{(pq)}\subset
\LQ_{\LA}.
\end{displaymath}
\end{re}

In the sequel we denote one of the two LC rings above by $R$, i.e.
$R =\LH$ or $\LQ_{pq}$. We also simply write $(\Delta, H, R, \beta, X(M,N))$
instead of $(\Delta_{R}, H_{R}, R_{R}, \beta_{R},X_{R}(M,N))$. Then
we have the following shocking result (a constructive mathematical
construction)
\begin{th}({\bf The existence of singular solutions of $AFE_{V}$}).

Let $V$ be a {\bf real vector space} with a {\bf continuous
idempotent endomorphism} $F : V \longrightarrow V$. Then {\bf each element
$v_{0} \in V$} is a {\bf $\pm$-fixed point of $F$}, i.e.
\begin{equation}
\;\;V\;=\;Fix(F),
\end{equation}
and moreover an arbitrary $v_{0} \in V$ has the following {\bf
Riesz-Bogoluboff-Kriloff Abstract Hodge Decomposition
(Representation)}($AHD_{RBK}$ in short)
\begin{equation}
(AHD_{RBK})\;v_{0}\;=\;\int \int_{\LS\times
X(M,N)}[(I \pm F)(v(\Delta^{2}(r),v_{0}))]d(\beta \otimes R)(s,r).
\end{equation}
\end{th}
{\bf Proof}. Since $R = \LH$ or ${\cal R} = \LQ_{pq}$  and according
to Lemma 2 , for each $v_{0}$ we have at our disposal the whole family
\begin{equation}
{\cal V}_{\pm}(v_{0})\;:=\;\{v_{\pm}(l,v_{0}) : l \ne \pm 1\}
\end{equation}
of solutions of the family of the {\bf abstract Fox equations}
\begin{equation}
v_{\pm}(l,v_{0})\;\pm \;l Fv_{\pm}(l,v_{0})\;=\;v_{0}.
\end{equation}
We substitute $l = \Delta^{2}(r), r \in  R^{*}, \Delta(r)\ne 1$,
in (3.81), thus obtaining
\begin{equation}
\frac{v_{\pm}(\Delta^{2}(r),v_{0})}{\Delta^{2}(r)}\;+\;F(v_{\pm}(\Delta^{2}(r),v_{0}))\;
=\;\frac{v_{0}}{\Delta^{2}(r)}.
\end{equation}
Integrating both sides of (3.82) with respect to the Haar measure $H$
on $X(M,N)$ and applying formula (3.67), here in the form
\begin{equation}
\int_{{\cal R}^{*}}f(r)\chi_{X(M,N)}(r)dH(r)\;=\;\int_{{\cal
R}^{*}}f(r^{-1})\frac{\chi_{X(M,N)}(r^{-1})}{\Delta^{2}(r)}dH(r),
\end{equation}
we obtain the equality
\begin{equation}
\int_{X(M,N)}\Delta^{-2}(r)v_{\pm}(\Delta^{2}(r),v_{0})dH(r)\;\pm \;
\int_{X(M,N)}\Delta^{-2}(r) F(v_{\pm}(\Delta^{-2}(r),v_{0}))dH(r)\;=\;
\end{equation}
\begin{displaymath}
=v_{0}\int_{X_{M,N}}\Delta^{-2}(r)dH(r) \;=:\;v_{0}m_{-2}(M,N),
\end{displaymath}
where we denote the $(-2)-R-{\bf moment}$ of the Haar measure $H$  on
the compact $X(M,N)$ by $m_{-2}(M,N)$.

Let us consider the expressions
\begin{equation}
\int_{X(M,N)}\frac{v_{\pm}(\Delta^{\pm2}(r),v)dH(r)}{\Delta^{2}(r)}.
\end{equation}
Applying the {\bf compact-$R$-Hilbert transform ${\cal H}$} in
the form of the Abstract Hodge Decomposition $(AHD_{R})$ :
\begin{equation}
\Delta(h^{-2})\;=\;\int_{\LS}\frac{dR(y)}{\Delta(h^{2}-y^{2})}\;=\;{\cal
 H}(R)(h^{2}),
\end{equation}
$h \in X(M,N)$, and using the {\bf Fubini theorem} we obtain
\begin{equation}
\int_{X(M,N)}\frac{v_{\pm}(\Delta^{\pm2}(r),v_{0})dH(r)}{\Delta(r^{2})}=
\int_{X(M,N)}v_{\pm}(\Delta^{\pm2}(r),v_{0})\int_{\LS}\frac{dR(y)}{\Delta(r^{2}-y^{2})}dH(r)=
\end{equation}
\begin{displaymath}
=\int_{\LS} dR(y)
\int_{X(M,N)}\frac{v_{\pm}(\Delta^{\pm2}(r),v_{0})dH(r)}{\Delta(r^{2}-y^{2})}.
\end{displaymath}
But, according to the formula $(W_{i})$, we can write the second inner integral
in the iterated integral above in the form : (since $\Delta(y)=1,
r/y=:r^{\prime}$)
\begin{equation}
\int_{X(M,N)}\frac{v_{\pm}(\Delta^{\pm2}(\frac{r}{y}),v_{0})dH(r)}
{\Delta^{2}(y)\Delta(1-(\frac{r}{y})^{2}}=
\Delta(y^{-3})\int_{X_{M,N}}\frac{v_{\pm}(\Delta^{\pm2}(r),v_{0})dH(r)}
{\Delta(1-r^{2})}\;,\;y \in \LS.
\end{equation}

But $\frac{dH(r)}{\Delta(1-r^{2})} =: dHer(r) = d\beta(r)$ is the {\bf Herbrand
distribution} of the invertion $I=I_{{\cal R}}$ of ${\cal R}^{*}$, i.e.
\begin{equation}
\int_{X(M,N)}\frac{v(\Delta^{2}(r),v_{0})dH(r))}{\Delta(1-r^{2})}\;=\;\int_{X(M,N)}
\frac{v(\Delta^{-2}(r),v_{0})dH(r)}{\Delta(1-r^{2})},
\end{equation}
since, for each integrable function $\phi$ the following is true
\begin{displaymath}
\int_{X(M,N)}\phi(I(r))dHer(r)\;=\;\int_{X(M,N)}\phi(r)d(I^{*}Her)(r)\;=
\end{displaymath}
\begin{displaymath}
\;=\;\int_{X(M,N)}\phi(r)dHer(r).
\end{displaymath}
Let us set:
\begin{equation}
v_{1}^{\pm}\;:=\;\int_{\LS}\Delta(y^{-3})dR(y)\int_{X(M,N)}\frac{v(\Delta^{2}(r)
,v_{0})dH(r)}{\Delta(1-r^{2})}.
\end{equation}
Since our "manipulations" only acted up on the parameters $l$  and under
the assumption, $F$ is continuous and linear, then we finally
obtain the RBK-integral representation above, which at the same time is
the {\bf singular solution} of $(AFE_{V})$ :
\begin{equation}
v_{1}^{\pm}\;\pm\;F(v_{1}^{\pm})\;=\;m_{-2}v_{0}.
\end{equation}

\begin{re}({\bf On a shocking consequence of the construction of Th.2. The
mathematics and logic}).

Obviously, the thesis (3.78) is not true (on the ground of classical
logic) for the majority of idempotent pairs $(V, F)$. Reely, let $V =
\LC$ be considered as the 2-dimensional Banach space over $\LR$ and let
$F = c$ be the complex conjugation. Then
\begin{displaymath}
\;\;Fix(c,\LC)\;=\;\LR \ne \LC.
\end{displaymath}

The construction in Th.2 is a following step in the old and well-known
philosophical problem : what is the connection between maths and
(classical) logic?

As it is well-known, {\bf Frege} saw mathematics as only a part of logic
(more exactly, according to Frege, the whole of mathematics can be reduced to
logic).

Probably the first mathematician, who questioned Frege's approach
to mathematics was pre-intuitionist {\bf Kronecker}. He attacked
well-known {\bf Cantor's proof} (in "naive" set theory), of the
existence of {\bf transcendental numbers} $t \in T$.

Let $\tilde{\LQ}$ be the (algebraically closed) field of {\bf algebraic
numbers} and assume that $TnD$ is true :
\begin{displaymath}
v_{Cl}(p \lor \sim p)\;=\;1 .
\end{displaymath}
We can write $TnD$ in quantifier form as the following true statement on
$\tilde{\LQ}$ (according to the laws of the quantifier  calculus) :
\begin{displaymath}
(C)\forall(x \in \LR)(x \in \tilde{\LQ}) \lor \sim(\forall (x \in
\LR)(x \in \tilde{\LQ}))\;=\;
\end{displaymath}
\begin{displaymath}
\;=\;\forall(x \in \LR)(x \in \tilde{\LQ})\lor (\exists(x \in \LR)(t
\in T)).
\end{displaymath}

Under the assumption, that the first term in the alternative $(C)$ is
true ($\tilde{\LQ}$ is countable!), it follows that $\LR$ should be countable,
which is impossible, according to the well-known Cantor theorem.

Thus, according to the rules of classical calculus of statements
and predicators, the second term of the alternative $(C)$ is true.
Thus, transcendental numbers exist.

But Cantor's reasoning does not give any information regarding a real
number, which is transcendental. In other words, it does not provide a
{\bf construction} of such a number.

According to Kronecker, the non-constructive character of Cantor's
proof of the existence of transcendental numbers is bad and hence its conclusion
should be rejected. But $(C)$ is only a specification of TnD. Thus,
questioning $(C)$ is identical to questioning TnD.
The immediate consequence of this was the rejection of classical
logic and construction of intuitione logic by {\bf Heyting}.
{\bf Brouwer} built constructive mathematics on this basis and showed that,
in general, many constructions violate TnD. For example Brouwer's construction
of the {\bf diagonal set of positive integers} $D\LN$ (the simplest {\bf Post
system} generated by a constructive object $\mid$ and the format
$\frac{x,x}{x}$ (cf.[ML, Sect.7])) violates the statement
\begin{displaymath}
(n \in D\LN)\lor (n \notin D\LN).
\end{displaymath}
Similarly, in our case the statement
\begin{displaymath}
(v_{0} \in Fix(F))\lor (v_{0}\notin Fix(F))
\end{displaymath}
violates TnD (a real infinity exists but no a potential infinity?)

The construction in Th.2 is an example of such a construction. In
reality, it leaveas out assumption : $v_{0} \in Fix(F)$. It seems that
it is much worse. {\bf It gives a contradiction in mathematics}.

According to {\bf Poincare}, the only thing, which we must demand
from an object which exists in mathematics is {\bf
non-contradictivity}. On the other hand, {\bf Godel's  well known result}
states that it is not possible to prove the { \bf non-contradictivity of
arithmetics of $\LN$} (and , in fact, the majority of axiomatic systems).
Moreover, the problem of the non-contradictivity of ZFC-set theory is much
more complicated than for such arithmetics. Thus (according to {\bf Gentzen's
 non-finistic proof} of the non-contradictivity of arithmetics), we can only
{\bf believe} that set theory is non-contradictory. But a belief is only
a belief, and for example, the proof of Th.2 seems be done properly, according
to classical logic, but it leads to classical mathematical
contradiction.

The only explanation of this phenomenon is the following : we use the methods
of {\bf measure theory} strictly, which is subsequently based on set theory,
in a strict manner. But according to the above discusion can this be ...
(contradictory)?

It is also very surprising, that such logical problems from the fundaments of
mathematics appeared during the work on the Riemann hypothesis. Maybe this is
one of the reasons that (RH) was unproven for so long and shows that (RH) is
not a standard mathematical problem.

Finally, all the logical problems with (RH) mentioned above should lead and
stimulate a subsequence discussion on mathemathical philosophy, very
similar to the discourses after {\bf Appel-Haken's proof} of the {\bf four
colour conjecture} (proved with help of a computer program). Can we accept a
proof of RH which is based - in its generality - on a theorem which leads to
a contradiction although, if we bound the domain of objects to  some "admissible"
$v_{0} \in Fix(V)$, then the construction is acceptable.
\end{re}

We now apply our theorem in the case $V = {\cal S}(\LR^{n})$ and $F =
{\cal F}_{n}$. Let $A^{+} = A_{n}^{+}(x)$  be a {\bf generalized amplitude},
i.e. any function from ${\cal S}(\LR^{n})$ with $A^{+}(0) = 0$. Then, according
to Th.2, there exists a {\bf RH-fixed point} $\omega_{A}^{+}$ (associated with
$A^{+}$, cf. [$M_{A}, Th.2$]), i.e.
\begin{equation}
(\omega^{+}_{A}\;-\;G)(x)\;=\;A^{+}(x)\;\;;\;x \in \LR^{n}.
\end{equation}
As we remarked in [$M_{A}$, Remark 15] (see also (3.91)), $\omega^{+}_{A}$
{\bf cannot exist} if $A$ is not a fixed point of ${\cal F}_{n}$, i.e.
according to Remark 3.

In [$M_{H}$, Prop.2]  we showed that a  direct solution of the {\bf RH-eigenvalue
problem} exists. We constructed a concrete example of the {\bf hermitian
amplitude} $A =A^{4}_{h_{0}}$  being an eigenvalue of the parameetrized Fourier
transform ${\cal F}_{h_{0}}$ and the {\bf RH-eigenvector} $\omega^{+}_{A}$ as the
{\bf fourth order hermite function} (cf.[$M_{H}$, (93.18)])
\begin{displaymath}
\omega^{+}_{A}(x)\;:=\;
H^{4}_{h_{0}}(x)\;:=\;h_{0}^{4}e^{-h_{0}^{2}x^{2}}
(16h_{0}^{4}x^{4}-48h_{0}^{2} x^{2}\;+\;12),
\end{displaymath}
which satisfies the equation
\begin{displaymath}
\omega^{+}_{A}(x)\;-\;12h_{0}^{4}e^{-h_{0}^{2}x^{2}}\;=:\;A(x)\;,\;x \in \LR.
\end{displaymath}
Here $h_{0} = \frac{\sqrt{3}}{2}$ is an amplitude parameter.
But it is very difficult ( either we cannot or it is not possible) to find
a direct analytic example of an {\bf RH-amplitude}. The main difficulty
is to find two fixed-points of the Fourier transform, which are both
stricly decreasing for $x >1$. In other words, the restriction : ${\cal
F}_{n}(A^{+}) = A^{+}$  is very restrictive.

For the parameter $p$ dependent Fourier transform ${\cal F}_{p}(f)(x)
:= \int_{\LR}e^{2p^{2}ixy}f(y)dy$, we showed in [$M_{H}$] that a direct
solution of the (-)RH-eigenvector problem exists.

Defining the (-)RH-eigenvector $\omega_{A}^{-}$ as the {\bf sixth order
hermitian function} (see $[M_{H}]$ for details)
\begin{displaymath}
\omega_{A}^{-}(x)\;:=\;H_{h_{0}}^{6}(x)\;=\;16h_{0}^{6}e^{-h_{0}^{2}x^{2}}(4
h_{0}^{6}x^{6}-30h_{0}^{4}x^{4}+45h_{0}^{2}x^{2}-7.5),
\end{displaymath}
we can define the amplitude $A^{-}$ by the formula :
\begin{displaymath}
A^{-}(x)\;:=\;\omega_{A}^{-}(x)\;+\;60h_{0}^{4}H^{2}_{h_{0}}(x).
\end{displaymath}

In the last part of this paper the fundamental role plays the {\bf
second canonical Hermite function}
\begin{displaymath}
H_{2}(x)\;=\;\pi G(x)(4 \pi x^{2}\;-\;1),\;\;x \in \LR.
\end{displaymath}
Integrating by parts twicely, we obtain that $H_{2}$ is a {\bf minus
fixed point} of ${\cal F}$ : $\hat{H}_{2}(x) \;=\; -H_{2}(x)$.

Then, according to Th.2, there exists a (-)RH-fixed point $\omega_{A}^{-}$
(associated with an amplitude $A^{-}$), i.e.
\begin{displaymath}
(\omega_{A}^{-}\;+\;H_{2})(x)\;=\;A^{-}(x)\;,\;x \in \LR.
\end{displaymath}

According to (3.91) $\omega_{A}^{-}$ cannot exists if $A^{-}$ is not a
minus fixed point of ${\cal F}_{1}$. But, if we take an amplitude $A^{-}$ in
such a way that the support of $(A^{-}\;-\;H_{2})$ :
\begin{displaymath}
supp(A^{-}\;-\;H_{2})\;=:\;S_{A}
\end{displaymath}
is the completion of a set with positive Lebesgue measure $\lambda_{n}$,
i.e. $\lambda_{n}(S_{A}^{c})>0$, then, according to the {\bf "separation of
variables"} construction from $Th.2$ we obtain
\begin{equation}
supp(\omega_{A}^{-})\;=\;S_{A}\;=\;and\;{\cal F}_{1}(\omega_{A}^{-})(x)\;=
\;\int_{S_{A}}cos(2\pi ixy)\omega_{A}^{-}(y)dy\;=:\;C_{A}(\omega_{A}^{-})(x).
\end{equation}
Since $C_{A}(H_{2}) \ne -H_{2}$, i.e. $H_{2}$ is not a {\bf minus -fixed point} of
$C_{A}$, then the calculation
\begin{displaymath}
C_{A}(A^{-})=C_{A}(\omega_{A}^{-}+H_{2})=C_{A}(\omega_{A}^{-})+C_{A}(H_{2})
\ne -(\omega_{A}^{-}+H_{2})= -A^{-}
\end{displaymath}
shows that , in this case, the notion of RH-fixed point does not lead to
a {\bf contradiction} and can exist for an amplitude $A^{-}$, which is not
the minus-fixed point of ${\cal F}_{1}$ (antinomies cannot be treated as a
threat to the fundaments of maths or logic).
\begin{re}
{\bf A. Wawrzy\,nczyk} (see [Wa, 3.8, Exercise 1d]) as well as we (see
[$M_{H}$, Prop.1] and [$M_{P}$, Remark 1]) have considered the following example
of the {\bf minus-fixed point} $K_{2}$ of ${\cal F}_{1}$ : in the
considerations concerning the {\bf quantum harmonic oscillator} in quantum
mechanics - one of the main roles plays the following {\bf second
Hermite function}
\begin{displaymath}
K_{2}(x)\;:=\;2\pi e^{-\pi x^{2}}(2 \pi x^{2}\;-\;1)
\end{displaymath}
with the property : $\hat{K}_{2}(x) = - K_{2}(x)$.

However, according to {\bf P. Biane}, since $G(x) := e^{-\pi x^{2}}$ is
a {\bf fixed point} of the canonical Fourier transform ${\cal F}_{1}$,
then integrating by parts twicely we obtain
\begin{displaymath}
\int_{\LR}e^{2\pi i xy}G^{\prime \prime}(y)dy\;=\;-2\pi ix
\int_{\LR}e^{2 \pi ixy}G^{\prime}dy\;=
\end{displaymath}
\begin{displaymath}
\;=\;-4\pi^{2}x^{2}\int_{\LR}e^{2\pi
ixy}G(y)dy\;=\;-(4\pi^{2}x^{2}G(x)\;-\;\pi G(x))-\pi G(x).
\end{displaymath}
Since
\begin{displaymath}
G^{\prime \prime}(x)\;=\;2\pi G(x)(2 \pi x^{2}\;-\;1)\;=\;K_{2}(x)
\end{displaymath}
then
\begin{displaymath}
\hat{H}_{2}(x):=\hat{(G^{\prime \prime})+\pi
G(x)}(x)\;=\;-4\pi^{2}x^{2}G(x)+\pi G(x)=-H_{2}(x),
\end{displaymath}
i.e. $H_{2}(x) := \pi G(x)(4 \pi x^{2}\;-\;1)$ is {\bf also} the {\bf
minus fixed point} of ${\cal F}_{1}$!

Since
\begin{displaymath}
H_{2}(x)(:=\omega_{A}^{-}(x))\;-\;K_{2}(x))\;=\;\pi G(x) (=: A^{-}).
\end{displaymath}

In the sequel, we call the equality : $H_{2}\;-\;K_{2} = \pi G$ - the
{\bf BMW-example}.

Since $A^{-}:=\pi G$ {\bf is not evidently the minus-fixed point of}
${\cal F}_{1}$ (since it is the +fixed-point of ${\cal F}_{1}$), then
the BMW-example :
(1). confirms the {\bf correctness} of our Th.2.-construction and
RH-fixed point paradox of Remark 3.
(2) We are not in possibility to explain that {\bf phenomena} on the
ground of the classical logic.!
\end{re}

\section{An (-)RH-fixed point proof of the generalized Riemann hypothesis}

As opossed to the purely algebraic notion of the quasi-fixed point $Q_{l}(v)$
considered in the previous section, here the main part is played by a
{\bf purely analytic} notion of the ($1$-dimensional) {\bf amplitude} $A$
(cf.[$M_{A}$] and [$M_{H}]$.

\begin{de}
We say that a function $A :\LR_{+}\longrightarrow \LR_{+}^{*}$ is an
{\bf PCID-amplitude}, if it is {\bf positive, continuous, integrable}
and {\bf (strictly) decreasing} on $[1,\infty)$.
\end{de}

The importance of PCID-amplitudes (amplitudes in short) follows from the fact
that an {\bf analytic Nakayama type lemma} holds for them (i.e. some very simple
analytic statement, trivial in proof - but powerful consequences). This lemma
establishes the sign of the action on the amplitude of the {\bf
plus-sine operator} $S_{+}$(cf.$[M_{A},Lemma 4]$ and $[M_{H}, Lemma2]$):
\begin{equation}
S_{+}(A)(a)\;:=\;\int_{0}^{\infty}A(x)sin(ax)dx \;,\;a>0.
\end{equation}
\begin{lem}
For each {\bf amplitude} $A$ and {\bf frequency} $a \in \LR^{*}_{+}$ the
following holds
\begin{equation}
\;\;S_{+}(A)(a)\;>\;0.
\end{equation}
\end{lem}

\begin{pr}({\bf On the positivity of the Rhfe(Ace)-trace} $Tr_{-}$).

Let $n = [k:\LQ]$ and  $A^{-}=A_{n}^{-}(x), x \in \LR^{n}$ from ${\cal S}
(\LR^{n})$ be a such function that for each $e$ from the fundamental domain
$E(k)$, the function $t \longrightarrow exp(t)A_{n}^{-}(exp(t)e) =: A_{n}^{e}(t)$
is a (1-dimensional) {\bf amplitude}.
Then, for each complex number $s = u +iv$ with $u \in (1/2, 1]$ and $v>0$,
the following {\bf Casteulnovo-Serre-Weil inequality} (CWS in short) holds
\begin{equation}
(CWS)\;Tr_{-}(\zeta_{k},A_{n}^{x})(s)\;:=\;\int_{1}^{+\infty}(t^{u}+t^{1-u})
sin(vlog t) \theta_{k}(A_{n}^{x})(t)dt\;>\;0.
\end{equation}
\end{pr}
{\bf Proof}. According to (2.35), (2.25) and (2.43) we have
\begin{equation}
Tr_{-}(\zeta_{k},A_{n})(s)=\sum_{0\ne I \in {\cal I}_{k}}\sum_{u \in
U(k)}\sum_{\xi \in R(I)}
\int_{1}^{\infty}(\int_{E(k)}A_{n}((N(I)^{-1}t)^{1/n} C(\xi u)\cdot e)
dH_{r}^{0}(e))
\end{equation}
\begin{displaymath}
(t^{u}+t^{1-u})sin(vlogt)dt/w(k).
\end{displaymath}
Let us denote the vector
\begin{displaymath}
C(I)x^{t}\;:=\;N(I)^{-1/n}[\sum_{j=1}^{n}x_{j}C_{1}(\alpha_{j}), ... ,
\sum_{j=1}^{n}x_{j}C_{n}(\alpha_{j})].
\end{displaymath}
After the substitution $t = e^{r}$ and changing of variables according
to the n-dimensional substitution : $e^{\prime} = e \cdot C(I)x^{t}$ we
obtain that
\begin{displaymath}
Tr_{-}(\zeta_{k}, A_{n})=\sum_{0\ne I \in {\cal I}_{k}}
\sum_{m \in \LZ^{n}}\frac{1}{w(k)\Delta_{r}(C(I)x^{t})}(\int_{E(k)}(
\int_{1}^{+\infty}A_{n}^{e}(e^{r/n})(e^{r(u+1)}+e^{r(2-u)})sin(vr)dr)
dH_{r}^{0}(e).
\end{displaymath}
Let us consider 1-dimensional {\bf amplitudes} of the form
\begin{displaymath}
{\cal A}^{e}_{n}(r):=e^{ru}(1+e^{r(1-2u)}) A_{n}^{e}(e^{r/n}).
\end{displaymath}
Since $\frac{d}{dr}(1+e^{r(1-2u)})<0$, if $u \in (1/2, 1]$, then under
our asumptions on the amplitude $A_{n}$ the function ${\cal A}^{e}_{n}(r)$ is
strictly decreasing. According to Lemma 4,
\begin{equation}
S_{+}({\cal A}^{e}_{n})(v)\;>\;0.
\end{equation}
Combining (4.97) with (4.98) we obtain the Proposition.

\begin{re}
The considered in Prop.3 the {\bf minus-trace} $Tr_{-}(\zeta_{k},A_{n}^{+})(s)$
is obviously associated with a {\bf minus-fixed} points $\omega^{-}$ of
${\cal F}_{n}$. We have seen that its {\bf positivity} is an
immediately consequence of the mentioned above analytic Nakayama lemma
(Fresnel lemma).

Instead of $Tr_{-}(\zeta_{k},A_{n}^{+})(s)$, in the first version of the
manuscript , we have considered the {\bf plus-trace} (associated with
$\omega^{+}$)
\begin{displaymath}
Tr_{+}(\zeta_{k},
A_{n}^{+})(s)\;:=\;\int_{1}^{\infty}(t^{u}-t^{1-u})sin(vlogt)\theta_{k}(A
_{n}^{+})(t)dt,
\end{displaymath}
(which obviously only differs from $Tr_{-}(\zeta_{k},A_{n}^{-})(s)=:Tr_{-}$ by
a {\bf sign} in the subintegral expression).

In opposite to the case of $Tr_{-}$ - the {\bf positivity} of $Tr_{+}>0$
- as it was independly communicated to the author by the private communications
by {\bf S. Albeverio}, {\bf P. Biane} and {\bf Z. Brze\,zniak}!, is not
an immediately consequence of the Fresnel lemma. In particular, the result
: $Tr_{+}>0$ needs a machine of stochastic analysis and is a final effect
of the existence of the so called {\bf Hodge measure} $H_{2}^{*}$ on
$\LC^{++}:=\{z \in \LC : Re(s)>0, Im(s)>0\}$, which gives the {\bf Laplace
representation} of the {\bf inverse of the Haar module of} $\LC$ :
\begin{displaymath}
\mid z \mid^{-2}\;=\;\int_{\LC^{++}}e^{z \cdot w}dH_{2}^{*}(w)\;\;,z \in
\LC^{++}.
\end{displaymath}
The existence of $H_{2}^{*}$ is far non-obvious. Even worse, many
peoples suggested to the author, that such the measure cannot exists!

Fortunately, the problem has a positive solution, although it is a very
technical and complicated in details result. So, we are not going to do
it in this paper.
\end{re}

\begin{lem}
$M(H_{2})(s)$ has no roots in the domain $\LC -\{1,2\}$ and has of order
1 poles at the points :2,1,0, -2,-4, ... .
\end{lem}
{\bf Proof}. Let us again recall that the canonical second Hermite function
$H_{2}(x)$ has the form (see also [$M_{H}$, Remark 1]) :
\begin{displaymath}
H_{2}(x)\;=\;\pi G(x)(4 \pi x^{2}\;-\;1).
\end{displaymath}
It is easy to check (integrating by parts, see [$M_{H}$, Prop. 7]) that
\begin{displaymath}
M(H_{2})(s)\;=\;(s-1)(s-2)M(G)(s-2)\;\;,\;for\;Re(s)>2.
\end{displaymath}
Since $M(H_{2})(s)$ is well-defined for $Re(s)>0$ (because $H_{2} \in
{\cal S}(\LR)$), then the above formula gives the analytic continuation
of the previous right-hand side formula - defined for $Re(s)>2$.

Making the substitution $\pi x^{2} =t$ in Gamma integral we
obtain
\begin{displaymath}
M(G)(s)=\int_{0}^{\infty}x^{s-2}e^{-\pi x^{2}}dx=\frac{\pi^{1-s}}
{2\pi}\Gamma(s/2),
\end{displaymath}
where $\Gamma$ denotes the classical gamma function. Since $M(G)(s)$
{\bf does not vanish anywhere}, then $M(H_{2})$ {\bf does not vanish}
for $\LC_{2} := \LC -\{1,2\}$.

The final preliminary result which is very convenient when we work with
$(gRH_{k})$ is the elegant {\bf Rouche theorem}(cf.e.g. [Ma, Th. XV.18 ])
: let $\Omega\subset \LC$ be a {\bf domain} and $D \subset \Omega$ be
{\bf compact}.
Let $f$ and $g$ be {\bf holomorphic functions} on $\Omega$, which
satisfy the two following {\bf Rouche's border conditions}:
\begin{equation}
f(z)\;\ne \;0\;for\;z \in \partial D\; (the\;border\;of\;D),
\end{equation}
and
\begin{equation}
\mid g(z) \mid\;<\;\mid f(z) \mid\;for\;z \in \partial D .
\end{equation}
Then the number of zeros $N_{D}(f+g)$ of the sum $f+g$ in $D$
(weighed by their orders) is equal to the number of zeros $N_{D}(f)$
of $f$ in $D$ ({\bf Rouche's thesis}, the adic type behaviour of
$N_{D}$), i.e.
\begin{equation}
N_{D}(f+g)\;=\;N_{D}(f).
\end{equation}
\begin{pr}({\bf A Rouche choice of the amplitude $A^{+}$ and lack of
roots of $\Gamma_{r}(G+A^{+})$}).
We can choose a plus amplitude $A_{n}^{+}$ in such a way that :
(1) the construction of the plus RH-fixed point $\omega_{A}^{+}$ in
Th.2 fulfills all the rigours of the classical logic, i.e. it does not
violate TnD.

(2) $\Gamma_{r}(G+A_{n}^{+})(s) \ne 0$ for $Re(s)>0$.
\end{pr}
{\bf Proof}. We use the Rouche theorem in the case : $\Omega = \LC$,

\begin{displaymath}
D\;=\;D_{M}\;:=\;\{s \in \LC: Re(s)\in [0,1], Im(s) \in [-M, M]\},\;M>0
\end{displaymath}
and
\begin{displaymath}
f(z)\;=\;\Gamma_{r}(G)(z)\;\;,\;\;g(z)\;=\;\Gamma_{r}(A_{n}^{+})(z).
\end{displaymath}
Since obviously $\Gamma_{r}(G)(z) \ne 0$  for $z \in \partial D_{M}$ and
$N_{D_{M}}(\Gamma_{r}(G)) = 0$, then it suffices to show that
\begin{equation}
\mid \Gamma_{r}(A_{n}^{+})(z) \mid^{2}\;<\;\mid \Gamma_{r}(G)(z) \mid^{2}
\end{equation}
for $z \in D_{M}$, to conclude that $N_{D_{M}}(G+A_{n}^{+}) = 0$  in $D_{M}$.

The inequality (4.102) is obviously equivalent to the inequality
\begin{equation}
Re^{2}(\int_{G_{r}}mod_{r}(g)^{s-1}A_{n}^{+}(g)d^{n}g)\;+\;Im^{2}(\int
_{G_{r}}mod_{r}(g)^{s-1}A_{n}^{+}(g)d^{n}g)<
\end{equation}
\begin{displaymath}
<\;Re^{2}(\int_{G_{r}}mod_{r}(g)^{s-1}G_{n}(g)d^{n}g)\;+\;Im^{2}(\int
_{G_{r}}mod_{r}(g)^{s-1}G_{n}(g)d^{n}g .
\end{displaymath}

Let us consider the Taylor expansion of $G_{n}(x) = e^{-\pi \mid \mid x
\mid \mid^{2}}$
\begin{displaymath}
G_{n}(x)\;=\;\sum_{m=0}^{\infty}\frac{(-1)^{m}\pi ^{2m}\mid \mid x \mid
\mid_{n}^{2m}}{m!},\;for\;x \in \LR^{n},
\end{displaymath}
and let us denote $g_{m} := \frac{\pi^{2m}}{m!}$.

Without loss of generality we can assume that $A_{n}^{+}$ is
NCID-amplitude, i.e. is {\bf negative} continuous integrable and such
that $-A_{n}^{+}$ is strictly decreasing for $\mid \mid x \mid
\mid_{n}\ge 1$. Reely, taking $s$ with $Im(s)<0$ we obtain :
$Tr_{+}(\zeta_{k}, -A_{n}^{+})>0$. We then can define $A_{n}^{+}$ as
follows
\begin{equation}
A_{n}^{+}(x)\;:=\;-G_{n}(x)\;\;for\;\;\mid \mid x \mid \mid_{n} \ge 1,
\end{equation}
and
\begin{equation}
A_{n}^{+}(x)\;:=\;-\sum_{m=2}^{\infty}(-1)^{m}g_{m-2} \mid\mid x
\mid\mid_{n}^{2m}\;\;if\;\;\mid\mid x \mid\mid_{n} \in [0,1].
\end{equation}

Since $\sum_{m=2}^{\infty}(-1)^{m}g_{m-2} = \sum_{m=0}^{\infty}g_{m}$,
then $A_{n}^{+}$ is {\bf continuous} and hence - NCID-amplitude.

Moreover from the definition we get that the support of $(G\;+\;A^{+})$
is the {\bf unit ball} $B_{n}$ of $\LR^{n}$.

Thus, to obtain (4.103) it suffices to show that for $Re(s) \ge 0$ holds
\begin{equation}
\mid Re(\int_{G_{r}\cap B_{n}}mod_{r}(g)^{s}A_{n}^{+}(g)dH_{r}(g))
\mid<\mid Re(\int_{G_{r}\cap B_{n}}mod_{r}^{s}G_{n}(g)dH_{r}(g)) \mid,
\end{equation}
and for $Im(s) \le 0$ holds
\begin{equation}
\mid Im(\int_{G_{r}\cap B_{n}}mod_{r}(g)^{s}A_{n}^{+}(g)dH_{r}(g)) \mid
<\mid Im(\int_{B_{n}\cap G_{r}}mod_{r}(g)^{s}G_{n}(g)dH_{r}(g)) \mid,
\end{equation}
since, according to the definition of $A_{n}^{+}$ we have
\begin{equation}
\int_{B_{n}^{c}\cap
G_{r}}mod_{r}(g)^{s}A_{n}^{+}(g)dH_{r}(g)=-\int_{B_{n}^{c}\cap
G_{r}}mod_{r}^{s}(g)G_{n}(g)dH_{r}(g)
\end{equation}
and - let us recall (see (2.16) and (2.17)) -
\begin{equation}
\Gamma_{r}(\mid\mid \cdot
\mid\mid_{n}^{2m}\chi_{B_{n}})(s)=\int_{G_{r}\cap
B_{n}}mod_{r}(g)^{s}\mid\mid g \mid\mid_{n}^{2m}dH_{r}(g)\;=
\end{equation}
\begin{displaymath}
\int_{G_{r}\cap B_{n}}mod_{r}(g)^{s}\mid\mid g
\mid\mid_{n}^{2m}\frac{d^{n}g}{mod_{r}(g)}=\int
\int_{(\LR_{+}^{*}\times G_{r}^{0})\cap
B_{n}}mod_{r}^{s}(t^{1/n}c)t^{2m}\mid\mid c
\mid\mid_{n}^{2m}\frac{d^{n}c dt}{mod_{r}(c)t}=
\end{displaymath}
\begin{displaymath}
\;\;=:\;\frac{c_{2m}}{s+2m},
\end{displaymath}
since
\begin{displaymath}
log(mod_{r}(g))=\sum_{i=1}^{r_{1}}log \mid x_{i}
\mid\;+\;\sum_{j=1}^{r_{2}}log \mid z_{j} \mid^{2} \le C(n) \mid\mid g
\mid\mid_{n}^{2}.
\end{displaymath}
It is obvious that for $Re(s)\ge 0$ we have
\begin{displaymath}
\sum_{m=2}^{\infty}\frac{(-1)^{m}c_{m}g_{m}(Re(s)+2m)}{\mid s+2m
\mid^{2}}<\sum_{m=0}^{\infty}\frac{(-1)^{m}c_{m}g_{m}(Re(s)+2m)}{\mid
s+2m \mid^{2}},
\end{displaymath}
since $1\;-\; \frac{\pi(x+2)}{\mid x+2 \mid^{2}}>0$, according to the
fact that the quadratic polynomial $x^{2}+(4-\pi)x+(4-\pi)>0$ for all
$x>0$. Analogously, for $Im(s)\le 0$ we have
\begin{displaymath}
-Im(s)\sum_{m=2}^{\infty}\frac{(-1)^{m}c_{m}g_{m}}{\mid s+2m \mid^{2}}<
-Im(s)\sum_{m=0}^{\infty}\frac{(-1)^{m}c_{m}g_{m}}{\mid s+2m \mid^{2}}.
\end{displaymath}

Thus, according to the definition of $A_{n}^{+}(G)$, from those strict
inequalities above, we claim (deduce) that the pair $(\Gamma_{r}(G_{n}),
\Gamma_{r}(A_{n}^{+}))$ satisfies the {\bf strong} Rouche boundary
conditions (4.99) and (4.100) on every compact set $D_{M}, M>0$ (and
not only on $\partial D_{M}$):
\begin{displaymath}
\mid \Gamma_{r}(A_{n}^{+})(s)\mid\;<\;\mid
\Gamma_{r}(G_{n})(s)\mid\;,\;s \in D_{M}.
\end{displaymath}
Converging with $M$ to the infinity we finally obtain
\begin{displaymath}
N_{D_{\infty}}(G_{n}\;+\;A_{n}^{+})\;=\;N_{D_{\infty}}(G_{n})\;=\;0.
\end{displaymath}

\begin{pr}({\bf A non-contradictory choice of the amplitude $A^{-}$ and
deleting of the problem of vanishing of $M(A^{-}-H_{2})$}).

We can choose an amplitude $A_{n}^{-}$ in such a way that :

(1) the construction of the (-)RH-fixed point $\omega_{A_{n}}^{-}$ in Th.2
fulfills all the rigours of classical logic, i.e. it does not violate
TnD.

(2) Even when $\Gamma_{r}(H_{2}-A^{-})(s)$ has {\bf zeros} in $Re(s)>0$
then still holds the (Face$_{-}$):
\begin{equation}
\Gamma_{r}(\omega_{A}^{-})(s)\zeta_{k}(s)\;=\;
\frac{\lambda_{k}}{s(s-1)}\;+\;
\end{equation}
\begin{displaymath}
\;+\;\int_{1}^{\infty}\int_{E}
\theta_{E}(\omega_{A}^{-})(ct^{1/n})(t^{s}\;+\;t^{1-s})dH_{r}^{0}(c)
\frac{dt}{t}(=:\int_{1}^{\infty}\Theta_{k}(\omega_{A}^{-})(t)(t^{s-1}+
t^{-s})dt).
\end{displaymath}
\end{pr}
{\bf Proof}. Let us consider the McLaurin expansion of $H_{2}^{n}$
\begin{displaymath}
H_{2}^{n}(x)\;=\;-\pi\;-\;\sum_{m=1}^{\infty}\frac{(-1)^{m}(-\pi)^{m+1}(4m+1)
\mid\mid x \mid\mid_{n}^{2m}}{m!},
\end{displaymath}
and let us denote $h_{m}:= \frac{\pi^{m+1}(4m+1)}{m!}$.

For a convenience of the reader, we give here all needed in the sequel
facts concerning the {\bf graph} of $H_{2}$ (it can be easy obtained by
using the elementary differential calculus). Thus :
\begin{equation}
H_{2}(0)=-\pi,\;H_{2}(\frac{1}{2\sqrt{\pi}})=0,\;H_{2}(1)=\pi
e^{-\pi}(4\pi-1)>0.
\end{equation}
Moreover, the function $H_{2}(x)$ is {\bf positive} for $x \ge 1/2\sqrt{\pi}$
and {\bf strictly decreasing} for $x \ge \sqrt{\frac{5}{2}}$. Finally,
the sequence $\{h_{m}\}$ is strictly decreasing for $m \ge 4$ (see also
[AM, Lemma 2]).

Looking at the graph of $H_{2}$ on $\LR_{+}$, we see that we can find
such $x_{2}>\sqrt{5/2}>1>x_{1}>1/2\sqrt{\pi}$ (since $H_{2}(x_{2})
\longrightarrow +\infty$ if $x_{2}\longrightarrow \infty$), that the
defined below function $A_{n}^{-}$ is an {\bf PCID-amplitude} :
\begin{displaymath}
A_{n}^{-}(x)\;:=\;H_{2}^{n}(x)\;\;if\;\; \mid\mid x \mid\mid_{n} \ge x_{2},
\end{displaymath}
and
\begin{displaymath}
A_{n}^{-}(x)\;:=\;L(x)\;if\;\mid\mid x \mid\mid_{n} \le x_{2},
\end{displaymath}
where by $L$ we denoted the line which connects the points $(x_{2},H_{2}
(x_{2}))$ and $(x_{1}, H_{2}(x_{1}))$  with $H_{2}(x_{1})>H_{2}(x_{1})$.
Moreover $(H_{2}\;-\;A_{n}^{-})(x)= 0$ for $\mid \mid x \mid \mid_{n}>x_{2}$.

The construction of an amplitude - let us say $A_{n}^{--}$ - with the property
that $\Gamma_{r}(H_{2}-A_{n}^{--})(s) \ne 0$ if $Re(s)>0, Im(s)>0$, i.e.
such $A_{n}^{--}$ that we could apply to it the Rouche theorem is much more
technically complicated (although possible). Therefore we are not going
to do it in this paper because we can overcome that problem as follows :
let us observe that Th.1 gives in fact a stronger result, i.e. it holds
{\bf without any assumption} on the vanishing of $\Gamma_{r}(\omega_{A}^{-})$.
Reely, beside the fact that we have not any exact information on the
zero-dimensional manifold $\Gamma_{r}(\omega_{A}^{-})(\LC):=\{s \in \LC :
\Gamma_{r}(A_{n}^{-})(s)=0\}$, the meromorphic functions :
$\Gamma_{r}(A_{n}^{-})(s)\zeta_{k}(s)$ and
$\int_{1}^{\infty}(t^{s-1}-t^{-s})\Theta_{k}(\omega_{A}^{-})(t)dt$ are
{\bf well-defined} for $Re(s)>0$ and - according to (Face) - {\bf
coincides} for $Re(s)>1$. Hence, according to the uniqueness of the
continuation of the analytic functions in regions - they must be equal
in $Re(s)>0$.

\begin{th}({\bf Existence of $gRhfe_{k}^{-}$}).

A pair of two ${\bf \Gamma \theta sinlog-factors} (F_{id},F_{c})$ indexed
by the Galois group $Gal(\LC/\LR) = \{id, c\}$ and another pair $(f_{1},
f_{2})$ of {\bf $\theta$sinlog-factors} satisfying
\begin{equation}
f_{1}(s)\;+\;f_{2}(s)\;\ne\;0\;for\;Re(s)\in (1/2,1]
\end{equation}
exist, such that the following $gRhfe_{k}^{-}$ ( with rational term $I$ and
the action of $Gal(\LC/\LR)$) holds
\begin{equation}
Im(\sum_{g \in
Gal(\LC/\LR)}(F_{g}\zeta_{k})(g(s))\;=\;\frac{\lambda_{k}(f_{1}(s)+f_{2}(s))}
{\mid s(s-1) \mid}I(s).
\end{equation}
\end{th}
{\bf Proof}. (I). {\bf The derivation of $gRhfe_{k}^{-}$}.

Let $a_{2}>a_{1}>0$ be arbitrary {\bf artificially chosen
$\zeta_{k}$-Cramer initial condition} and let $s=u+iv=Re(s)+iIm(s)$ be
fixed. We consider a simple {\bf non-homogeneous system} of two linear
equations in two variables $p_{1}$ and $p_{2}$ of the form :
\begin{equation}
p_{1}v(u-1)\;+\;p_{2}vu\;=\;a_{1}\;-\;a_{2}
\end{equation}
\begin{displaymath}
p_{1}vu\;+\;p_{2}v(u-1)\;=\;a_{2}\;-\;a_{1}.
\end{displaymath}
This system is a {\bf Cramer system}, iff $s$ does not belong to
the algebraic $\LR$-variety $I(\LC)$. The main determinant of (4.104) is
$I(s)$ and its solution is given by the formulas
\begin{equation}
p_{1}\;=\;p_{1}(Im(s))\;=\;\frac{(a_{2}\;-\;a_{1})}{v}\;>\;0
\end{equation}
and
\begin{equation}
p_{2}\;=\;p_{2}(Im(s)) \;=\;\frac{(a_{1}\;-\;a_{2})}{v}\;=-p_{1}<\;0.
\end{equation}
Let $A^{-}$ be an {\bf amplitude} chosen according to the Proposition 5. Then
according to Proposition 5, there exists a {\bf (-)RH-fixed point}
$\omega_{A}^{-}$ of ${\cal F}_{n}$, i.e.
\begin{equation}
\omega^{-}_{A}\;+\;H_{2}\;=\;A^{-}.
\end{equation}
In the sequel we simply write $\omega_{1} = \omega^{-}_{A}$.
We denote the standard n-dimensional second Hermite (-)fixed point of
${\cal F}_{n}$ by $\omega_{2} = H_{2} = H_{2}^{n}$.

We set (cf.(2.48))
\begin{equation}
J_{i}(s)\;:=\;\int_{1}^{+\infty}(t^{u-1}\;+\;t^{-u})sin(vlogt)\Theta_{k}
 (\omega_{i})(t)dt\;,\;i=1,2.
\end{equation}

The integrals $J_{i}$ above are {\bf quasi-invariant} under the substitutions
: $t = x^{r}, r>0$, i.e. the substitution $t = x^{p_{1}v}, v>0$
gives
\begin{equation}
J_{1}(s)=p_{1}v
\int_{1}^{\infty}(x^{p_{1}v(u-1)}+x^{-p_{1}vu})sin(p_{1}v^{2}logx)\Theta
_{k}(\omega_{1})(x^{p_{1}v})x^{(p_{1}v-1)}dx=:J_{1}^{r}(s)
\end{equation}
In the same way, the substitution $t=x^{-p_{2}v}, v>0$ gives
\begin{equation}
J_{2}(s)=-p_{2}v\int_{1}^{\infty}(x^{-p_{2}v(u-1)}+x^{p_{2}vu})sin(-p_{2}
v^{2}logx)\Theta_{k}(\omega_{2})(x^{-p_{2}v})x^{-(p_{2}v+1)}dx=:
J_{2}^{r}(s).
\end{equation}
Thus, the equalities $J_{i}(s) = J_{i}^{r}(s), i=1,2$ hold on the
domain $\{s \in \LC: Im(s) \ge 0\}$. But obviously the integrals are
imaginary parts of the {\bf analytic} function $\Gamma_{(r_{1},r_{2})}(\omega
_{i})\zeta_{k}-\lambda_{k}/W$ defined on $\LC-\{0,1\}$. Hence, they
must be equal everywhere. In particular, the second equality is {\bf
invariant} to the operation of {\bf complex conjugation} $c$, i.e.
\begin{equation}
J_{2}(c(s))=p_{2}v\int_{1}^{\infty}(x^{p_{2}v(u-1)}+x^{p_{2}vu})sin(-p_{2}v^{2}
log x)\Theta_{k}(\omega_{2})(x^{p_{2}v})x^{p_{2}v-1}dx=J_{2}^{r}(c(s)).
\end{equation}

Since $\omega_{i} \in {\cal S}(\LR^{n})$, for each $q>1$ we have
\begin{displaymath}
max_{x\ge 1}\mid x^{q} \Theta_{k}(\omega_{i})(x^{p_{i}v})\mid
<\infty.
\end{displaymath}

According to the {\bf elementary mean value theorem}, there exists such
an $x_{i} = x_{i}(s,a_{1},a_{2}) \in [1,\infty)$ and $q = q(a_{1},a_{2},u)>1$
that
\begin{equation}
J_{i}(c_{i}(s))=p_{i}vsin((-1)^{i+1}p_{i}v^{2}logx_{i})x^{q}_{i}\Theta_{k}
(\omega_{i})(x_{i}^{p_{i}v})\int_{1}^{\infty}(x^{p_{i}v(u-1)-a_{i}}+
x^{-p_{i}vu-a_{i}})x^{-q}dx
\end{equation}
\begin{displaymath}
\;=:\;f_{i}(s)\int_{i}(s),
\end{displaymath}
where $c_{1}=id$ and $c_{2}=c$.

The number $q$ is obviously chossen in such a way that the integrals
$\int_{i}(s)$ are convergent.

Using the (Face) (cf.(2.36)) and the nation from (4.122) we obtain
\begin{equation}
Im((\Gamma_{(r_{1},r_{2})}(\omega_{i})\zeta_{k})(c_{i}(s))=
\frac{\lambda_{k}I(c_{i}(s))}{\mid s(s-1)
\mid^{2}}+f_{i}(s)\int_{i}(s),
\end{equation}
or equivalently
\begin{equation}
Im((\Gamma_{(r_{1},r_{2})}(\omega_{1})f_{2}\zeta_{k}))(s)=\frac{(f_{2}I)(s)}
{\mid s(s-1) \mid^{2}}\;+\;(f_{1}f_{2})(Im(s))\int_{1}(s),
\end{equation}
together with
\begin{equation}
Im(\Gamma_{(r_{1},r_{2})}(\omega_{2})f_{1}\zeta_{k})(c(s))=\frac{-(f_{1}I)
(s)}{\mid s(s-1) \mid^{2}}\;+\;(f_{1}f_{2})(Im(s))\int_{2}(s).
\end{equation}

By defining the {\bf $\Gamma \theta$ sinlog-factors} as
\begin{equation}
F_{id}(s):=(\Gamma_{(r_{1},r_{2})}(\omega_{1}f_{2}))(s)\;and\;F_{c}(s):=
(\Gamma_{(r_{1},r_{2})}(\omega_{2}f_{1}))(s),
\end{equation}
and substrating (4.125) from (4.124), according to the choice of the pair
$(p_{1},p_{2})$ in (4.114) (which is the solution of the Cramer system) we
finally obtain $(gRhfe_{k}^{-})$.

(II).{\bf Positivity of $Tr(\zeta_{k},A)$} (It is a very subtle "game"
of signs - on the bourder of subtlety) .

According to the construction of $\omega^{A}$,  we have
\begin{equation}
A\;=\;\omega_{1}\;+\;\omega_{2}.
\end{equation}
By Proposition 3 on the positivity of the trace, we have
\begin{equation}
0<Tr_{-}(\zeta_{k},A)(s)\;=\;(J_{1}+J_{2})(s)\;=\;J_{1}(s)+J_{2}(m(s))\;=
\end{equation}
\begin{displaymath}
J_{1}(s)\;+\;J_{2}(c(a(s))),
\end{displaymath}
where - for a moment - we denoted the affinic antyconjugation as
$a(s):= (1-u)+iv$, and
\begin{equation}
J_{2}(s)\;=\;Im(\int_{1}^{\infty}(t^{s-1}-t^{-s}))\Theta_{k}(t)dt= -J_{2}
(m(s)).
\end{equation}
Moreover, on the basis of the notation in (4.122) we have
\begin{displaymath}
f_{2}(s)\;=\;\frac{J_{2}(c(s))}{\int_{2}(s)},
\end{displaymath}
and therefore
\begin{equation}
f_{2}(a(s))\;=\;-f_{2}(s)\;and\;\int_{2}(a(s))\;=\;\int_{2}(s).
\end{equation}

Since the pair $(p_{1},p_{2})$ is the solution of the Cramer system
(4.104), we obtain
\begin{equation}
-\int_{2}(s)\;=\;\int_{1}(s).
\end{equation}
Hence, combining (4.128), (4.130) and (4.131) we finally obtain
\begin{equation}
0<Tr_{-}(\zeta_{k},A)(s)=f_{1}(s)\int_{1}(s)\;+\;f_{2}(a(s))\int_{2}(a(s))\;
=
\end{equation}
\begin{displaymath}
\;\;\;=\;\int_{1}(s)(f_{1}(s)\;+\;f_{2}(s)),
\end{displaymath}
i.e.
\begin{displaymath}
\;\;\;f_{1}(s)\;+\;f_{2}(s)\;\ne\;0\;for\;Re(s)\in (1/2,1]
\end{displaymath}
which proves Theorem 3.
\begin{re}
It is a very exciting fact that to prove ($gRH_{k}$) we need only {\bf
two} functional equations for $\zeta_{k}(s)$!, whereas - among number
theory specialists - we have met with the quite opposite opinion - that
even infinitely many f.e. for $\zeta_{\LQ}(s)$ are {\bf not sufficient}
to proof (RH)! (e.g. H. Iwaniec).
\end{re}

Obviously $(gRhfe_{k}^{-})$ immediately implies the {\bf generalized
Riemann Hypothesis}. Assume that there exists a zero $s_{0}$ of
$\zeta_{k}$ in the set $\{s \in \LC: Re(s \in (1/2,1]\}$. Then
\begin{displaymath}
\sum_{g \in Gal(\LC/\LR)}(F_{g}\zeta_{k})(g(s_{0}))\;=\;0,
\end{displaymath}
since, according to HRace, the zeros of zeta lie symmetrically with
respect to the lines : $Im(s) = 0$ and $Re(s) = 1/2$. But, on the other
hand, we have
\begin{displaymath}
\frac{(f_{1}\;+\;f_{2})(s_{0})}{\mid s_{0}(s_{0}-1)
\mid^{2}}I(s_{0})\;\ne\;0,
\end{displaymath}
which is impossible according to $(gRhfe_{k}^{-})$.

\begin{re}
The CWS-inequality
\begin{displaymath}
Tr_{Gal(\LC/\LR)}^{k}(s)\;=\;Tr^{k}_{G}(s)\;:=\;\frac{\lambda_{k}(f_{1}+
f_{2})(s)}{\mid s(s-1) \mid^{2}}\;>\;0,
\end{displaymath}
is {\bf exceptional}(fundamental) to the proof of $(gRH_{k})$. That is
very surprising that similar kinds of positivity conditions (explored
also in [$M_{A}$], $[M_{H}]$ and [AM]) are strictly connected with (RH) :

In $[M_{CG}]$, based on $[M_{L}]$ we showed that the positivity of the Cauchy-Gaussian trace
$Tr_{CG}$ implies the Riemann hypothesis.

In [B] de Branges showed that the positivity of his trace $Tr_{B}$
would imply the Riemann hypothesis (also in the case of some
$L$-functions).

Below we briefly remind the reader that the positivity of the {\bf Weil trace}
$Tr_{W}$ leads to the Riemann hypothesis.

As it is well-known (cf.e.g.[L, XVII.3]), A. Weil formulated an
equivalent form of the Riemann hypothesis (the {\bf Weil Formula} (WF
in short)) in terms of the {\bf positivity} of his functional: let
${\cal SB}(\LR)$ be the {\bf restricted Barner-Schwartz space} of all
functions of the form
\begin{displaymath}
F(x)\;=\;P(x)e^{-Kx^{2}}
\end{displaymath}
with some real constant $K>0$ and some polynomial $P$ (cf.[L,XVII.3]).
Then ${\cal SB}(\LR)$ is self-dual, and functions from this space
satisfy the {\bf three Barner conditions} (cf.[L]) : finitness of
variation, Dirichlet normalization and asymptotic symmetry at zero.

For each $F \in {\cal SB}(\LR)$ is well-defined its {\bf conjugation}
\begin{displaymath}
F^{*}(x)\;:=\;F(-x),
\end{displaymath}
and $F$ is of {\bf positive type} if $F$ is equal to its {\bf Rosatti
convolution}
\begin{displaymath}
F\;=\;F_{0} * F_{0}^{*},
\end{displaymath}
for some $F_{0} \in {\cal SB}(\LR)$. (So we see that ${\cal SB}(\LR)$
is also closed under the convolution $*$).
For $s = \sigma +it$, we can consider the {\bf two-sided
Laplace-Fourier transform}
\begin{displaymath}
\hat{F}(s)\;:=\;\int_{\LR}F(x)e^{(1/2-\sigma)x}e^{itx}dx,
\end{displaymath}
and the {\bf Weil functional} $W$ defined as
\begin{displaymath}
W_{k}(\Phi)\;:=\;\sum_{\rho, \zeta_{k}(\rho)=0, Im(\rho)\ne 0}\Phi(\rho).
\end{displaymath}
Then the Riemann hypothesis is equivalent to the positivity of Weil's trace
\begin{equation}
Tr_{W}(F_{0})\;:=\;W_{k}(F_{0}* F^{*}_{0})\;\;\ge\;0,
\end{equation}
for all $F_{0} \in {\cal SB}(\LR)$.

{\bf Weil's condition} is much more general.

Let $k$ be a number field, $\chi$ a {\bf Hecke character}, ${\cal f}_{\chi}$
the {\bf conductor}, ${\cal D}$ the {\bf local different} and $d_{\chi} =
N({\cal D}f_{\chi})$.

Let us consider the $L^{*}_{k}$-function
\begin{equation}
L_{k}^{*}(s;\chi)\;:=\;[(2\pi)^{-n(k)}2^{r_{1}}d_{\chi}]^{s/2}\prod_{v
\in S_{\infty}(k)}\Gamma(s_{v}/2)L(s;\chi),
\end{equation}
where $L(s;\chi)$ is the {\bf Hecke $L$-function} associated with
$\chi$, i.e. the usual product over unramified prime ideals for $\chi$
and $s_{v} := N_{v}(s+i\phi_{v})+ \mid m_{v} \mid$ (cf.[L]).
{\bf Weil's functional} $W$ in this case is obviously the sum
\begin{displaymath}
W_{L}(F)\;=\;\sum_{L(\rho, \chi)=0}F(\rho)\;\;,\;\;F\in {\cal SB}(\LR).
\end{displaymath}
In short, the {\bf generalized Riemann hypothesis} for $L_{k}^{*}(\cdot;\chi)$,
$gRH_{k}(\chi)$, states that $Re(\rho)=\frac{1}{2}$ for all zeros $\rho$ of
$L_{k}(\cdot;\chi)$ in the critical strip.
Well-known {\bf Weil's theorem}(cf.[L,Th.3.3]) asserts that $gRH_{k}(\chi)$ is
equivalent to the property that
\begin{equation}
(WC)\forall(F_{0} \in {\cal SB}(\LR))(W_{L}(F_{0}*F_{0})\ge 0).
\end{equation}
\end{re}

In particular, we have thus proved Weil's theorem for the Dedekind
zetas.
\begin{th}
For all $F_{0}$ in the restricted Schwartz space ${\cal SB}(\LR)$ the
following holds
\begin{displaymath}
W_{\zeta_{k}}(F_{0}* F_{0}^{*})\;\ge \;0.
\end{displaymath}
\end{th}

\section{The generalized Riemann hypothesis for all Dirichlet
L-functions}

A first generalization of the Riemann zeta function comes from
{\bf Dirichlet}[Di], who for a character $\chi$ of $(\LZ/m \LZ)^{*}$,
that is, a homomorphism from $(\LZ/ m \LZ)^{*}$ to $\LC^{*}$, considered
the series
\begin{equation}
L(s, \chi)\;:=\;\sum_{n=1}^{\infty}\frac{\chi(n)}{n^{s}},
\end{equation}
where $\chi(n):= \chi([n])$ for $(n,m)=1$ and $\chi(n)=0$ for $(n,m)\ne
1$. He used these L-series to prove his theorem on primes in arithmetic
progressions, in which of principal importance is the fact that the
value of $L(s, \chi)$ is nonzero at the point $s=1$.

Let $m$ be a natural number amd $\zeta_{m}$ a primitive mth root of
unity, that is, a complex number with $\zeta_{m}^{m}=1$ and
$\zeta_{m}^{i}\ne 1$ for $1 \le i \le m$. In this section we consider
extensions $k$ that arise from $\LQ$ through the adjunction of roots of
unity. The field $k = \LQ(\zeta_{m})$ is called the {\bf mth cyclotomic
field}, since as points in the complex plane they divide the circle
into equal arcs (see [K, Sect. 6.4]).

Since the development by Kummer of the theory of cyclotomic fields (see
e.g. [K]) one proves $L(1, \chi)\ne 0$ for characters $\chi$ different
from the {\bf trivial character} $\chi_{0}$ ( $L(s, \chi_{0})$ has a simple
pole at $s=1$) most naturally with the help of the following result
(see [K, Sect.8.2, Th.8.2.1.]) :

 for any integer $ m \in \LN$
\begin{equation}
\zeta_{\LQ(\zeta_{m})}(s)\;=\;\prod_{p \mid
m}(1\;-\;\frac{1}{N(p)^{s}})^{-1}\prod_{\chi}L(s, \chi),
\end{equation}
where the right-hand product runs over all characters of $(\LZ/ m
\LZ)^{*}$.

\begin{th}({\bf $gRH_{m}$ for Dirichlet L-functions})
Let $m$ be any positive integer and $\chi_{m} : \LF_{m}^{*}= (\LZ/ m\LZ)^{*}
\longrightarrow  \LC$ any character of the multiplicative group of the
finite ring $\LF_{m}$. Let also $\chi_{m}$ be corresponding {\bf
Dirichlet character}. Then the following implication is true :
\begin{equation}
(gRH_{m})\;If\;L(s, \chi_{m})=0\;with\;Im(s)\ne
0\;then\;Re(s)=1/2.
\end{equation}
In particular, the {\bf Weil trace}  $Tr_{W,m}(F_{0}) :=
\sum_{\rho, L(\rho,\chi_{m})=0, Im(\rho)\ne 0 } (F_{0}(\rho)*F_{0}(\rho))$
associated with the L-function $L(s, \chi_{m})$ is {\bf positive}, i.e.
\begin{equation}
Tr_{W,m}(F_{0})\;\ge\;0,
\end{equation}
for all $F_{0}$ from the Barner-Schwartz space ${\cal SB}(\LR)$.
\end{th}
{\bf Proof}. Assume (a contrary) that there is a {\bf zero} $s_{0}$ of
$L(s, \chi_{m})$ in the domain : $Re(s)\in (0,1)-\{1/2\},
Im(s)\ne 0$ of $\LC$. Then, according to the "spliting formula" (5.137) we
obtain that
\begin{displaymath}
\zeta_{\LQ(\zeta_{m})}(s_{0})\;=\;0,
\end{displaymath}
what obviously is not possible according to $gRH_{k}$.

Thus, the generalized Riemann hypothesis for Dirichlet L-functions -
according to (5.137) - is directly and immediately reduced (or is the
consequence) of the generalized Riemann hypothesis for Dedekind zetas -
proved in the previous Section.

e-mail: madrecki@im.pwr.wroc.pl

\end{document}